\documentclass[draft]{amsart}

\usepackage{amsthm}
\usepackage[leqno]{amsmath}
\usepackage{latexsym,amsfonts,amssymb}
\usepackage[all]{xy} \SelectTips{eu}{} \SilentMatrices
\usepackage[urlcolor=blue]{hyperref}

\newcommand{\numberseries}{\mdseries}   

\newlength{\thmtopspace}                
\newlength{\thmbotspace}                
\newlength{\thmheadspace}               
\newlength{\thmindent}                  

\renewcommand{\subparagraph}{\vspace{\thmbotspace}}

\newtheoremstyle{bfupright head,slanted body}
                {\thmtopspace}{\thmbotspace}
                {\slshape}{\thmindent}{\bfseries}{.}{\thmheadspace}
                {{\numberseries \thmnumber{(#2) }}\thmnote{#3}}

\newtheoremstyle{bfupright head,upright body}
                {\thmtopspace}{\thmbotspace}
                {\upshape}{\thmindent}{\bfseries}{.}{\thmheadspace}
                {{\numberseries \thmnumber{(#2) }}\thmnote{#3}}

\newtheoremstyle{bfit head,upright body}
                {\thmtopspace}{\thmbotspace}
                {\upshape}{\thmindent}{\upshape}{.}{\thmheadspace}
                {{\numberseries\thmnumber{(#2) }}
                {\bfseries\itshape\thmnote{\negthickspace#3}}}

\newtheoremstyle{it head,upright body}
                {\thmtopspace}{\thmbotspace}
                {\upshape}{\thmindent}{\upshape}{.}{\thmheadspace}
                {{\numberseries\thmnumber{(#2) }}
                {\itshape\thmnote{\negthickspace#3}}}

\newtheoremstyle{fixed bf head,slanted body}
                {\thmtopspace}{\thmbotspace}{\slshape}
                {\thmindent}{\bfseries}{.}{\thmheadspace}
                {{\numberseries \thmnumber{(#2) }}\thmname{#1}\thmnote{ (#3)}}

\newtheoremstyle{fixed bf head,upright body}
                {\thmtopspace}{\thmbotspace}{\upshape}
                {\thmindent}{\bfseries}{.}{\thmheadspace}
                {{\numberseries \thmnumber{(#2) }}\thmname{#1}\thmnote{ (#3)}}

\newtheoremstyle{fixed bfit head,upright body}
                {\thmtopspace}{\thmbotspace}{\upshape}
                {\thmindent}{\bfseries\itshape}{.}{\thmheadspace}
                {{\numberseries \thmnumber{(#2) }}\thmname{#1}\thmnote{ (#3)}}

\newtheoremstyle{sc head,small body}
                {\thmtopspace}{\thmbotspace}
                {\small\upshape}{\thmindent}{\scshape}{.}{\thmheadspace}
                {\thmname{#1}}

\newtheoremstyle{numbered paragraph}
                {\thmtopspace}{\thmbotspace}{\upshape}
                {\thmindent}{\upshape}{}{0pt}
                {{\numberseries \thmnumber{(#2) }}}

\newtheoremstyle{unnumbered paragraph}
                {\thmtopspace}{\thmbotspace}{\upshape}
                {\parindent}{\upshape}{}{0pt}

\theoremstyle{bfupright head,slanted body}
\newtheorem{res}{}[section]             \newtheorem*{res*}{}

\theoremstyle{bfit head,upright body}
                 \newtheorem*{com*}{}

\theoremstyle{bfupright head,upright body}
\newtheorem{bfhpg}[res]{}               \newtheorem*{bfhpg*}{}

\theoremstyle{it head,upright body}
               \newtheorem*{ithpg*}{}

\theoremstyle{sc head,small body}

\theoremstyle{fixed bf head,slanted body}
\newtheorem{thm}[res]{Theorem}          \newtheorem*{thm*}{Theorem}
\newtheorem{prp}[res]{Proposition}      \newtheorem*{prp*}{Proposition}
\newtheorem{cor}[res]{Corollary}        \newtheorem*{cor*}{Corollary}
\newtheorem{lem}[res]{Lemma}            \newtheorem*{lem*}{Lemma}

\theoremstyle{fixed bf head,upright body}
\newtheorem{dfn}[res]{Definition}       \newtheorem*{dfn*}{Definition}
\newtheorem{con}[res]{Construction}     \newtheorem*{con*}{Construction}
      \newtheorem*{obs*}{Observation}
\newtheorem{rmk}[res]{Remark}           \newtheorem*{rmk*}{Remark}
\newtheorem{exa}[res]{Example}          \newtheorem*{exa*}{Example}
         \newtheorem*{exe*}{Exercise}
            \newtheorem{stp*}{Setup}

\theoremstyle{numbered paragraph}

\theoremstyle{unnumbered paragraph}
\newtheorem{ipg*}{}

\newlength{\thmlistleft}        
\newlength{\thmlistright}       
\newlength{\thmlistpartopsep}   
\newlength{\thmlisttopsep}      
\newlength{\thmlistparsep}      
\newlength{\thmlistitemsep}     

\newcounter{eqc}
\newenvironment{eqc}{\begin{list}{\upshape (\textit{\roman{eqc}})}%
    {\usecounter{eqc}%
      \setlength{\leftmargin}{\thmlistleft}%
      \setlength{\labelwidth}{\thmlistleft}%
      \setlength{\rightmargin}{\thmlistright}%
      \setlength{\partopsep}{\thmlistpartopsep}%
      \setlength{\topsep}{\thmlisttopsep}%
      \setlength{\parsep}{\thmlistparsep}%
      \setlength{\itemsep}{\thmlistitemsep}}}%
  {\end{list}}%

\newcounter{prt}
\newenvironment{prt}{\begin{list}{\upshape (\alph{prt})}%
    {\usecounter{prt}%
      \setlength{\leftmargin}{\thmlistleft}%
      \setlength{\labelwidth}{\thmlistleft}%
      \setlength{\rightmargin}{\thmlistright}%
      \setlength{\partopsep}{\thmlistpartopsep}%
      \setlength{\topsep}{\thmlisttopsep}%
      \setlength{\parsep}{\thmlistparsep}%
      \setlength{\itemsep}{\thmlistitemsep}}}%
  {\end{list}}%

\newcommand{\prtlbl}[1]{{\upshape(#1)}}

\newcounter{rqm}
  {\end{list}}%

  {\end{list}}%

  {\end{list}}%

\newcounter{exercise}
  {\end{list}}%

\newenvironment{prf*}[1][Proof]{%
  \begin{proof}[\bf #1]
    \setcounter{equation}{0}
    \renewcommand{\theequation}{\arabic{equation}}}
  {\end{proof}
}

\newcommand{\pgref}[1]{(\ref{#1})}

\newcommand{\thmref}[2][Theorem~]{#1\pgref{thm:#2}}
\newcommand{\corref}[2][Corollary~]{#1\pgref{cor:#2}}
\newcommand{\prpref}[2][Proposition~]{#1\pgref{prp:#2}}
\newcommand{\lemref}[2][Lemma~]{#1\pgref{lem:#2}}

\newcommand{\conref}[2][Construction~]{#1\pgref{con:#2}}

\newcommand{\exaref}[2][Example~]{#1\pgref{exa:#2}}
\newcommand{\rmkref}[2][Remark~]{#1\pgref{rmk:#2}}
\newcommand{\secref}[2][Section~]{#1\ref{sec:#2}}

\newcommand{\partpgref}[2]{(\ref{#1})\prtlbl{#2}}

\newcommand{\partcorref}[3][Corollary~]{#1\partpgref{cor:#2}{#3}}
\newcommand{\partprpref}[3][Proposition~]{#1\partpgref{prp:#2}{#3}}
\newcommand{\partlemref}[3][Lemma~]{#1\partpgref{lem:#2}{#3}}

\renewcommand{\eqref}[1]{\pgref{eq:#1}}

\newcommand{\thmcite}[2][?]{\cite[thm.~#1]{#2}}

\newcommand{\prpcite}[2][?]{\cite[prop.~#1]{#2}}
\newcommand{\lemcite}[2][?]{\cite[lem.~#1]{#2}}

\newcommand{\seccite}[2][?]{\cite[sec.~#1]{#2}}

\numberwithin{equation}{res}
\newcounter{marcom}

\newcommand{\lgtR}{\operatorname{length}R}
\newcommand{\bet}[3][R]{\beta^{#1}_{#2}(#3)}
\newcommand{\Rm}{(R,\m)}
\newcommand{\Rmk}{(R,\m,\k)}

\newcommand{\set}[2][\mspace{1mu}]{\{#1 #2 #1\}}
\newcommand{\setof}[3][\,]{\{#1#2 \mid #3#1\}}
\renewcommand{\Im}[1]{\nobreak{\operatorname{Im}#1}}
\newcommand{\Coker}[1]{\nobreak{\operatorname{Coker}#1}}
\newcommand{\Ker}[1]{\nobreak{\operatorname{Ker}#1}}
\newcommand{\lgt}[2][R]{\operatorname{length}_{#1}#2}
\newcommand{\dif}[2][]{{\partial}^{#2}_{#1}}
\newcommand{\qtext}[1]{\quad\text{#1}\quad}
\newcommand{\qqtext}[1]{\qquad\text{#1}\qquad}
\newcommand{\qand}{\qtext{and}}
\newcommand{\qqand}{\qqtext{and}}
\newcommand{\NN}{\mathbb{N}}
\newcommand{\ZZ}{\mathbb{Z}}
\newcommand{\QQ}{\mathbb{Q}}
\newcommand{\is}{\cong}
\renewcommand{\le}{\leqslant}
\renewcommand{\ge}{\geqslant}
\newcommand{\onto}{\twoheadrightarrow}
\newcommand{\lra}{\longrightarrow}
\newcommand{\xra}[2][]{\xrightarrow[#1]{\;#2\;}}
\renewcommand{\dim}[2][R]{\operatorname{dim}_{#1}#2}
\newcommand{\poly}[2][\k]{#1[#2]}

\newcommand{\m}{\mathfrak{m}}
\newcommand{\n}{\mathfrak{n}}
\newcommand{\p}{\mathfrak{p}}
\newcommand{\Hom}[3][R]{\operatorname{Hom}_{#1}(#2,#3)}
\newcommand{\Ext}[4][R]{\operatorname{Ext}_{#1}^{#2}(#3,#4)}
\newcommand{\tp}[3][R]{\nobreak{#2\otimes_{#1}#3}}

\newcommand{\Fc}[2]{L^{#1}(#2)}

\renewcommand{\k}{\mathsf{k}}
\newcommand{\Po}[2][i]{P^{(#1,#2)}}
\newcommand{\GL}[2][\k]{\operatorname{GL}_{#2}(#1)}

\newcommand{\subeqref}[2]{\eqref{#1}$_\mathrm{#2}$}
\newcommand{\mnotmm}{\m\backslash\m^2}
\newcommand{\ann}[1]{( 0 : #1)}
\newcommand{\card}[1]{\operatorname{card}(#1)}
\newcommand{\embdim}{\operatorname{emb.\!dim}}
\newcommand{\w}{\tilde{w}}
\newcommand{\x}{\tilde{x}}
\newcommand{\y}{\tilde{y}}
\newcommand{\z}{\tilde{z}}
\newcommand{\Alpha}{\textrm{A}}
\newcommand{\Beta}{\textrm{B}}
\newcommand{\Epsilon}{\textrm{E}}
\newcommand{\Eta}{\textrm{H}}
\newcommand{\Kappa}{\textrm{K}}
\newcommand{\Iota}{\textrm{I}}
\newcommand{\Chi}{\textrm{X}}
\newcommand{\J}{\textrm{J}}
\newcommand{\Zeta}{\textrm{Z}}

\def\urltilda{\kern -.15em\lower .7ex\hbox{\~{}}\kern .04em}
\def\widebardisplay#1{%
  \setbox0=\hbox{$\displaystyle #1$}
  \dimen0=\wd0%
  \advance\dimen0 by -3.2pt
  \vbox{%
    \nointerlineskip%
    \moveright 1.2pt 
    \vbox{\hrule width \dimen0}%
    \nointerlineskip%
    \kern 1.25pt
    \box0%
    }%
  }

\def\widebartext#1{%
  \setbox0=\hbox{$#1$}
  \dimen0=\wd0%
  \advance\dimen0 by -3.2pt
  \vbox{%
    \nointerlineskip%
    \moveright 1.2pt 
    \vbox{\hrule width \dimen0}%
    \nointerlineskip%
    \kern 1.25pt
    \box0%
    }%
  }

\def\widebarscript#1{%
  \setbox0=\hbox{$\scriptstyle #1$}
  \dimen0=\wd0%
  \advance\dimen0 by -2pt
  \vbox{%
    \nointerlineskip%
    \moveright 1pt 
    \vbox{\hrule width \dimen0}%
    \nointerlineskip%
    \kern .8pt
    \box0%
    }%
  }

\def\widebarscriptscript#1{%
  \setbox0=\hbox{$\scriptscriptstyle #1$}
  \dimen0=\wd0%
  \advance\dimen0 by -1pt
  \vbox{%
    \nointerlineskip%
    \moveright .5pt 
    \vbox{\hrule width \dimen0}%
    \nointerlineskip%
    \kern .6pt
    \box0%
    }%
  }

\def\widebar#1{\mathchoice%
  {\widebardisplay{#1}}%
  {\widebartext{#1}}%
  {\widebarscript{#1}}%
  {\widebarscriptscript{#1}}%
  }

\begin{document}

\allowdisplaybreaks[2]

\title{Brauer--Thrall for totally reflexive modules}

\author[L.\,W. Christensen]{Lars Winther Christensen}

\address{L.W.C.\newline\hspace*{1em} Department of Math.\ and Stat.,
  Texas Tech University, Lubbock, TX 79409, U.S.A.}

\email{lars.w.christensen@ttu.edu}

\urladdr{http://www.math.ttu.edu/\urltilda lchriste}

\author[D.\,A. Jorgensen]{David A. Jorgensen}

\address{D.A.J.\newline\hspace*{1em} Department of Mathematics,
  University of Texas, Arlington, TX~76019, U.S.A.}

\email{djorgens@uta.edu}

\urladdr{http://dreadnought.uta.edu/\urltilda dave}

\author[H. Rahmati]{Hamidreza Rahmati}

\address{H.R.\newline\hspace*{1em} Department of Math.\ and Stat.,
  Texas Tech University, Lubbock, TX 79409, U.S.A.}

\curraddr{Mathematics Department, Syracuse University, Syracuse, NY
  13244, U.S.A.}

\email{hrahmati@syr.edu}

\author[J. Striuli]{Janet Striuli}

\address{J.S.\newline\hspace*{1em} Department of Math.\ and C.S.,
  Fairfield University, Fairfield, CT~06824,~U.S.A.}

\email{jstriuli@fairfield.edu}

\urladdr{http://www.faculty.fairfield.edu/jstriuli}

\thanks{This research was partly supported by NSA grant H98230-10-0197
  (D.A.J.), NSF grant DMS\,0901427 (J.S.), and by a UNL Faculty
  Development Fellowship (R.W.)}

\author[R. Wiegand]{Roger Wiegand}

\address{R.W.\newline\hspace*{1em} Department of Mathematics,
  University of Nebraska, Lincoln, NE~68588, U.S.A.}

\email{rwiegand1@math.unl.edu}

\urladdr{http://www.math.unl.edu/\urltilda rwiegand1}

\date{\today}




\begin{abstract}
  Let $R$ be a commutative noetherian local ring that is not
  Gorenstein. It is known that the category of totally reflexive
  modules over $R$ is representation infinite, provided that it
  contains a non-free module. The main goal of this paper is to
  understand how complex the category of totally reflexive modules can
  be in this situation.

  Local rings $\Rm$ with $\m^3=0$ are commonly regarded as the
  structurally simplest rings to admit diverse categorical and
  homological characteristics. For such rings we obtain conclusive
  results about the category of totally reflexive modules, modeled on
  the Brauer--Thrall conjectures. Starting from a non-free cyclic
  totally reflexive module, we construct a family of indecomposable
  totally reflexive $R$-modules that contains, for every $n\in\NN$, a
  module that is minimally generated by $n$ elements. Moreover, if the
  residue field $R/\m$ is algebraically closed, then we construct for
  every $n\in\NN$ an infinite family of indecomposable and pairwise
  non-isomorphic totally reflexive $R$-modules, each of which is
  minimally generated by $n$ elements. The modules in both families
  have periodic minimal free resolutions of period at most $2$.
\end{abstract}

\maketitle

\tableofcontents

\section{Introduction and synopsis of the main results} 
\label{sec:intro}

\noindent
The representation theoretic properties of a local ring bear pertinent
information about its singularity type. A notable illustration of this
tenet is due to Herzog \cite{JHr78} and to Buchweitz, Greuel, and
Schreyer \cite{BGS-87}. They show that a complete local Gorenstein
algebra is a simple hypersurface singularity if its category of
maximal Cohen--Macaulay modules is representation finite.  A module
category is called \emph{representation finite} if it comprises only
finitely many indecomposable modules up to isomorphism. Typical
examples of maximal Cohen--Macaulay modules over a Cohen--Macaulay
local ring are high syzygies of finitely generated modules.

Over a Gorenstein local ring, all maximal Cohen--Macaulay modules
arise as high syzygies, but over an arbitrary Cohen--Macaulay local
ring they may not. Totally reflexive modules are infinite syzygies
with special duality properties; the precise definition is given
below. One reason to study these modules---in fact, the one discovered
most recently---is that they afford a characterization of simple
hypersurface singularities among all complete local algebras, i.e.\
without any \emph{a priori} assumption of Gorensteinness.  This
extension of the result from \cite{BGS-87,JHr78} is obtained in
\cite{CPST-08}. It is consonant with the intuition that the structure
of high syzygies is shaped predominantly by the ring, and the same
intuition guides this work.

The key result in \cite{CPST-08} asserts that if a local ring is not
Gorenstein and the category of totally reflexive modules contains a
non-free module, then it is representation infinite.  The main goal of
this paper is to determine how complex the category of totally
reflexive modules is when it is representation infinite. Our results
suggest that it is often quite complex; \thmref[Theorems~]{firstbt}
and \thmref[]{secondbt} below are modeled on the Brauer--Thrall
conjectures.

For a finite dimensional algebra $A$, the first Brauer--Thrall
conjecture asserts that if the category of $A$-modules of finite
length is representation infinite, then there exist indecomposable
$A$-modules of arbitrarily large length. The second conjecture asserts
that if the underlying field is infinite, and there exist
indecomposable\linebreak $A$-modules of arbitrarily large length, then
there exist infinitely many integers $d$ such that there are
infinitely many indecomposable $A$-modules of length $d$. The first
conjecture was proved by Ro\u\i ter~(1968); the second conjecture has
been verified, for example, for algebras over algebraically closed
fields by Bautista~(1985) and Bongartz~(1985).

\begin{center}
  $*\:*\:*$
\end{center}

\noindent
In this section, $R$ is a commutative noetherian local ring.  The
central questions addressed in the paper are: Assuming that the
category of totally reflexive $R$-modules is representation infinite
and given a non-free totally reflexive $R$-module, how does one
construct an infinite family of pairwise non-isomorphic totally
reflexive $R$-modules? And, can one control the size of the modules in
the family in accordance with the Brauer--Thrall conjectures?

A finitely generated $R$-module $M$ is called \emph{totally reflexive}
if there exists an infinite sequence of finitely generated free
$R$-modules
\begin{equation*}
  F\colon\quad \cdots \lra F_{1} \lra F_{0} \lra F_{-1} \lra \cdots,
\end{equation*}
such that $M$ is isomorphic to the module $\Coker{(F_1 \to F_0)}$, and
such that both $F$ and the dual sequence $\Hom{F}{R}$ are exact. These
modules were first studied by Auslander and Bridger \cite{MAsMBr69},
who proved that $R$ is Gorenstein if and only if every $R$-module has
a totally reflexive syzygy. Over a Gorenstein ring, the totally
reflexive modules are precisely the maximal Cohen--Macaulay modules,
and these have been studied extensively. In the rest of this section
we assume that $R$ is not Gorenstein.

Every syzygy of an indecomposable totally reflexive $R$-module is,
itself, indecomposable and totally reflexive; a proof of this folklore
result is included in \secref{1}. Thus if one were given a totally
reflexive module whose minimal free resolution is non-periodic, then
the syzygies would form the desired infinite family; though one cannot
exercise any control over the size of the modules in the family.

In practice, however, the totally reflexive modules that one typically
spots have periodic free resolutions.  To illustrate this point,
consider the $\QQ$-algebra
\begin{equation*}
  A=\poly[\QQ]{s,t,u,v}/(s^2,t^2,u^2,v^2,uv, 2su+tu, sv+tv).
\end{equation*}
It has some easily recognizable totally reflexive modules---$A/(s)$
and $A/(s+u)$ for example---whose minimal free resolutions are
periodic of period at most $2$. It also has indecomposable totally
reflexive modules with non-periodic free resolutions. However, such
modules are significantly harder to recognize. In fact, when Gasharov
and Peeva~\cite{VNGIVP90} did so, it allowed them to disprove a
conjecture of Eisenbud.

The algebra $A$ above has Hilbert series $1 + 4\tau + 3\tau^2$. In
particular, $A$ is a local ring, and the third power of its maximal
ideal is zero; informally we refer to such rings as \emph{short.}  For
these rings, \cite{DAJLMS06} gives a quantitative measure of how
challenging it can be to recognize totally reflexive modules with
non-periodic resolutions.

Short local rings are the structurally simplest rings that accommodate
a wide range of homological behavior, and \cite{VNGIVP90} and
\cite{DAJLMS06} are but two affirmations that such rings are excellent
grounds for investigating homological questions in local algebra.


\begin{center}
  $*\:*\:*$
\end{center}


\noindent
For the rest of this section, assume that $R$ is short and let $\m$ be
the maximal ideal of $R$. Note that $\m^3$ is zero and set $e =
\dim[R/\m]{\m/\m^2}$.  A reader so inclined is welcome to think of a
standard graded algebra with Hilbert series $1 + e\tau + h_2\tau^2$.

The families of totally reflexive modules constructed in this paper
start from cyclic ones. Over a short local ring, such modules are
generated by elements with cyclic annihilators; Henriques and
\c{S}ega~\cite{IBHLMS} call these elements \emph{exact zero divisors.}
The ubiquity of exact zero divisors in short local algebras is a
long-standing empirical fact. Its theoretical underpinnings are found
in works of Conca~\cite{ACn00} and Hochster and
Laksov~\cite{MHcDLk87}; we extend them in \secref{quad}.

The main results of this paper are
\thmref[Theorems~]{firstbt}--\thmref[]{secondbt}.  It is known from
work of Yoshino~\cite{YYs03} that the length of a totally reflexive
$R$-module is a multiple~of~$e$. In \secref{m3} we prove the existence
of indecomposable totally reflexive $R$-modules of every possible
length:

\begin{thm}[Brauer--Thrall I]
  \label{thm:firstbt}
  If there is an exact zero divisor in $R$, then there exists a family
  $\set{M_n}_{n\in\NN}$ of indecomposable totally reflexive
  $R$-modules with $\lgt{M_n} = ne$ for every $n$. Moreover, the
  minimal free resolution of each module $M_n$ is periodic of period
  at most $2$.
\end{thm}
\noindent
Our proof of this result is constructive in the sense that we exhibit
presentation matrices for the modules $M_n$; they are all upper
triangular square matrices with exact zero divisors on the diagonal.
Yet, the strong converse contained in \thmref{supplies} came as a
surprise to us. It illustrates the point that the structure of the
ring is revealed in high syzygies.

\begin{thm}
  \label{thm:supplies}
  If there exists a totally reflexive $R$-module without free
  summands, which is presented by a matrix that has a column or a row
  with only one non-zero entry, then that entry is an exact zero
  divisor in $R$.
\end{thm}

\noindent
These two results---the latter of which is distilled from
\thmref{ezdexist}---show that existence of totally reflexive modules
of any size is related to the existence of exact zero divisors. One
direction, however, is not unconditional, and in \secref{example} we
show that non-free totally reflexive $R$-modules may also exist in the
absence of exact zero divisors in $R$.  If this phenomenon appears
peculiar, some consolation may be found in the next theorem, which is
proved in \secref{quad}. For algebraically closed fields, it can be
deduced from results in \cite{ACn00,MHcDLk87,YYs03}.

\begin{thm}
  \label{thm:generic}
  Let $\k$ be an infinite field, and let $R$ be a generic standard
  graded $\k$-algebra which (1) has Hilbert series $1 + h_1\tau +
  h_2\tau^2$, (2) is not Gorenstein, and (3) admits a non-free totally
  reflexive module. Then $R$ has an exact zero divisor.
\end{thm}

If the residue field $R/\m$ is infinite, and there is an exact zero
divisor in $R$, then there are infinitely many different ones; this is
made precise in \thmref{BT21}. Together with a couple of other results
from \secref{BT2} this theorem yields:

\begin{thm}[Brauer--Thrall II]
  \label{thm:secondbt}
  If there is an exact zero divisor in $R$, and the residue field $\k
  =R/\m$ is algebraically closed, then there exists for each $n\in\NN$
  a family $\set{M_n^\lambda}_{\lambda \in \k}$ of indecomposable and
  pairwise non-isomorphic totally reflexive $R$-modules with
  $\lgt{M_n^\lambda}=ne$ for every $\lambda$. Moreover, the minimal
  free resolution of each module $M_n^\lambda$ is periodic of period
  at most $2$.
\end{thm}

\begin{center}
  $*\:*\:*$
\end{center}

\noindent
The families of modules in \thmref[Theorems~]{firstbt} and
\thmref[]{secondbt} come from a construction that can provide such
families over a general local ring, contingent on the existence of
minimal generators of the maximal ideal with certain relations among
them.  This construction is presented in \secref{1} and analyzed in
\secref[Sections~]{bt} and \secref[]{bt2}. To establish the
Brauer--Thrall theorems, we prove in \secref[Sections~]{m3} and
\secref[]{BT2} that the necessary elements and relations are available
in short local rings with exact zero~divisors. The existence of exact
zero divisors is addressed in \secref[Sections~]{exist},
\secref[]{quad}, and \secref[]{example}.  In \secref{ext} we give
another construction of infinite families of totally reflexive
modules. It applies to certain rings of positive dimension, and it
does not depend on the existence of exact zero~divisors.

\section{Totally acyclic complexes and exact zero divisors} 
\label{sec:1}

\noindent
In this paper, $R$ denotes a commutative noetherian ring.  Complexes
of $R$-modules are graded homologically. A complex
\begin{equation*}
  F:\quad \cdots \lra F_{i+1} \xra{\dif[i+1]{F}} F_{i}
  \xra{\dif[i]{F}} F_{i-1} \lra \cdots
\end{equation*}
of finitely generated free $R$-modules is called \emph{acyclic} if
every cycle is a boundary; that is, the equality $\Ker{\dif[i]{F}} =
\Im{\dif[i+1]{F}}$ holds for every $i\in\ZZ$. If both $F$ and the dual
complex $\Hom{F}{R}$ are acyclic, then $F$ is called \emph{totally
  acyclic}. Thus an $R$-module is totally reflexive if and only if it
is the cokernel of a differential in a totally acyclic complex.

\subparagraph The annihilator of an ideal $\mathfrak{a}$ in $R$ is
written $\ann{\mathfrak{a}}$. For principal ideals $\mathfrak{a} =
(a)$ we use the simplified notation $\ann{a}$.

Recall from \cite{IBHLMS} the notion of an exact zero divisor: a
non-invertible element $x \ne 0$ in $R$ is called an \emph{exact zero
  divisor} if one of the following equivalent conditions holds.
\begin{eqc}
\item There is an isomorphism $\ann{x} \is R/(x)$.
\item There exists an element $w$ in $R$ such that $\ann{x} = (w)$ and
  $\ann{w} = (x)$.
\item There exists an element $w$ in $R$ such that
  \begin{equation}
    \label{eq:epzdtac}
    \cdots \lra R \xra{w} R \xra{x} R \xra{w} R \lra \cdots
  \end{equation}
  is an acyclic complex.
\end{eqc}

\noindent
For every element $w$ as above, one says that $w$ and $x$ form an
\emph{exact pair} of zero divisors in $R$. If $R$ is local, then $w$
is unique up to multiplication by a unit in $R$.

\begin{rmk}
  \label{rmk:ezd}
  For a non-unit $x\ne 0$ the conditions $(i)$--$(iii)$ above are
  equivalent~to
  \begin{eqc}\setcounter{eqc}{3}
  \item There exist elements $w$ and $y$ in $R$ such that the sequence
    \begin{equation*}
      R \xra{w} R \xra{x} R \xra{y} R
    \end{equation*}
    is exact.
  \end{eqc}
  Indeed, exactness of this sequence implies that there are equalities
  $\ann{y} = (x)$ and $\ann{x} = (w)$. Thus there is an obvious
  inclusion $(y) \subseteq (w)$ and, therefore, an inclusion $\ann{w}
  \subseteq \ann{y} = (x)$. As $x$ annihilates $w$, this forces the
  equality $(x)=\ann{w}$. Thus $(iv)$ implies $(ii)$ and, clearly,
  $(iv)$ follows from $(iii)$.
\end{rmk}

Starting from the next section, we shall assume that $R$ is local. We
use the notation $\Rm$ to fix $\m$ as the unique maximal ideal of $R$
and the notation $\Rmk$ to also fix the residue field $\k= R/\m$.

A complex $F$ of free modules over a local ring $\Rm$ is called
\emph{minimal} if one has $\Im{\dif[i]{F}} \subseteq \m F_{i-1}$ for
all $i\in\ZZ$. Let $M$ be a finitely generated $R$-module and let $F$
be a minimal free resolution of $M$; it is unique up to
isomorphism. The $i$th \emph{Betti number}, $\bet{i}{M}$, is the rank
of the free module $F_i$, and the $i$th \emph{syzygy} of $M$ is the
module $\Coker{\dif[i+1]{F}}$.

\begin{rmk}
  \label{rmk:ezdtac}
  The complex \eqref{epzdtac} is isomorphic to its own dual, so if it
  is acyclic, then it is totally acyclic. Thus if $w$ and $x$ form an
  exact pair of zero divisors in $R$, then the modules $(w) \is R/(x)$
  and $(x) \is R/(w)$ are totally reflexive. Moreover, it follows from
  condition $(iv)$ that a totally acyclic complex of free modules, in
  which four consecutive modules have rank $1$, has the form
  \eqref{epzdtac}. This means, in particular, that over a local ring
  $R$, any totally reflexive module $M$ with $\bet{i}{M} = 1$ for $0
  \le i \le 3$ has the form $M \is R/(x)$, where $x$ is an exact zero
  divisor. For modules over short local rings we prove a stronger
  statement in \thmref{ezdexist}.
\end{rmk}

\enlargethispage*{\baselineskip}
The next lemma is folklore; it is proved for Gorenstein rings in
\cite{JHr78}.

\begin{lem}
  \label{lem:indec}
  Let $R$ be local and let $T$ be a minimal totally acyclic complex of
  finitely generated free $R$-modules. The following conditions are
  equivalent:
  \begin{eqc}
  \item The module $\Coker{\dif[i]{T}}$ is indecomposable for some
    $i\in\ZZ$.
  \item The module $\Coker{\dif[i]{T}}$ is indecomposable for every
    $i\in\ZZ$.
  \end{eqc}
  In particular, every syzygy of an indecomposable totally reflexive
  $R$-module is indecomposable and totally reflexive.
\end{lem}

\begin{prf*}
  For every $i\in\ZZ$ set $M_i = \Coker{\dif[i]{T}}$. By definition,
  the modules $M_i$ and $\Hom{M_i}{R}$ are totally reflexive and one
  has $M_i \is \Hom{\Hom{M_i}{R}}{R}$ for every $i\in\ZZ$. Assume that
  for some integer $j$, the module $M_j$ is indecomposable. For every
  $i<j$, the module $M_j$ is a syzygy of $M_i$, and every summand of
  $M_i$ has infinite projective dimension, as $T$ is minimal and
  totally acyclic. Thus $M_i$ is indecomposable. Further, if there
  were a non-trivial decomposition $M_i \is K \oplus N$ for some
  $i>j$, then the dual module $\Hom{M_i}{R}$ would decompose as
  $\Hom{K}{R} \oplus \Hom{N}{R}$, where both summands would be
  non-zero $R$-modules of infinite projective dimension. However,
  $\Hom{M_j}{R}$ is a syzygy of $\Hom{M_i}{R}$, so $\Hom{M_j}{R}$
  would then have a non-trivial decomposition and so would $M_j \is
  \Hom{\Hom{M_j}{R}}{R}$; a contradiction.
\end{prf*}

\begin{lem}
  \label{lem:tac}
  Let $T$ and $F$ be complexes of finitely generated free
  $R$-modules. If $T$ is totally acyclic, and the modules $F_i$ are
  zero for $i \ll 0$, then the complex $\tp{T}{F}$ is totally acyclic.
\end{lem}

\begin{prf*}
  The complex $\tp{T}{F}$ is acyclic by
  \lemcite[2.13]{CFH-06}. Adjointness of Hom and tensor yields the
  isomorphism $\Hom{\tp{T}{F}}{R} \is \Hom{F}{\Hom{T}{R}}$. As the
  complex $\Hom{T}{R}$ is acyclic, so is $\Hom{F}{\Hom{T}{R}}$ by
  \lemcite[2.4]{CFH-06}.
\end{prf*}

\begin{dfn}
  \label{dfn:TR}
  For $n \in \NN$ and elements $y$ and $z$ in $R$ with $yz=0$, let
  $\Fc{n}{y,z}$ be the complex defined as follows
  \begin{equation*}
    \label{eq:Fn}
    \Fc{n}{y,z}_i =
    \begin{cases}
      R\quad&\text{for $0 \ge i \ge -n+1$}\\
      0 &\text{elsewhere}
    \end{cases}
    \qand
    \dif[i]{\Fc{n}{y,z}} =
    \begin{cases}
      y\quad&\text{for $i$ even}\\
      -z &\text{for $i$ odd.}
    \end{cases}
  \end{equation*}
\end{dfn}

If $w$ and $x$ form an exact pair of zero divisors in $R$, and $T$ is
the complex \eqref{epzdtac}, then the differentials of the complex
$\tp{T}{\Fc{n}{y,z}}$ have a particularly simple form; see
\rmkref{del1}. In \secref[Sections~]{bt} and \secref[]{bt2} we study
modules whose presentation matrices have this form.

\begin{con}
  \label{con:indec}
  Let $n \in \NN$, let $\Iota_n$ be the $n \times n$ identity matrix,
  and let $r_i$ denote its $i$th row. Consider the $n \times n$
  matrices $\Iota^o_n$, $\Iota^e_n$, $\J^o_n$, and $\J^e_n$ defined by
  specifying their rows as follows
  \begin{align*}
    (\Iota^o_n)_i &=
    \begin{cases}
      r_i & \text{$i$ odd}\\
      0 & \text{$i$ even}
    \end{cases}
    & &\mspace{-20mu}\text{and}\mspace{-20mu}& (\Iota^e_n)_i &=
    \begin{cases}
      0& \text{$i$ odd}\\
      r_{i}& \text{$i$ even}
    \end{cases}\\
    (\J^o_n)_i &=
    \begin{cases}
      r_{i+1} & \text{$i$ odd}\\
      0 & \text{$i$ even}
    \end{cases}
    & &\mspace{-20mu}\text{and}\mspace{-20mu}& (\J^e_n)_i &=
    \begin{cases}
      0 & \text{$i$ odd}\\
      r_{i+1}& \text{$i$ even}
    \end{cases}
  \end{align*}
  with the convention $r_{n+1} =0$. The equality $\Iota^o_n+\Iota^e_n
  = \Iota_n$ is clear, and the matrix $\J_n=\J^o_n+\J^e_n$ is the $n
  \times n$ nilpotent Jordan block with eigenvalue zero.

  For elements $w$, $x$, $y$, and $z$ in $R$ let $M_n(w,x,y,z)$ be the
  $R$-module with presentation matrix $\Theta_n(w,x,y,z) =
  w\Iota^o_n+x\Iota^e_n+ y\J^o_n+z\J^e_n$; it is an upper triangular
  $n\times n$ matrix with $w$ and $x$ alternating on the diagonal, and
  with $y$ and $z$ alternating on the superdiagonal:
  \begin{equation*}
    \Theta_n(w,x,y,z)
    = \begin{pmatrix}
      w & y & 0 & 0 & 0 &  \dots\\[.3ex]
      0 & x & z & 0 & 0 &  \dots\\[.3ex]
      0 & 0 & w & y & 0 &  \dots\\[.3ex]
      0 & 0 & 0 & x & z &  \\[.3ex]
      0 & 0 & 0 & 0 & w &  \smash{\ddots}\\[-.6ex]
      \vdots & \vdots & \vdots & \vdots & & \smash{\ddots}
    \end{pmatrix}.
  \end{equation*}
  For $n=1$ the matrix has only one entry, namely $w$. For $n=2$, the
  matrix does not depend on $z$, so we set $M_2(w,x,y) = M_2(w,x,y,z)$
  for every $z\in R$.

  Note that if $\Rm$ is local and $w$, $x$, $y$, and $z$ are non-zero
  elements in $\m$, then $M_n(w,x,y,z)$ is minimally generated by $n$
  elements.
\end{con}

\begin{rmk}
  \label{rmk:del1}
  Assume that $w$ and $x$ form an exact pair of zero divisors in $R$,
  and let $T$ be the complex \eqref{epzdtac}, positioned such that
  $\dif[1]{T}$ is multiplication by $w$. Let $n\in\NN$ and let $y$ and
  $z$ be elements in $R$ that satisfy $yz=0$. It follows from
  \lemref{tac} that the complex $\tp{T}{\Fc{n}{y,z}}$ is totally
  acyclic. It is elementary to verify that the differential
  $\dif[i]{\tp{T}{\Fc{n}{y,z}}}$ is given by the matrix
  $\Theta_n(w,x,y,z)$ for $i$ odd and by $\Theta_n(x,w,-y,-z)$ for $i$
  even, cf.~\conref{indec}.  In particular, the module $M_n(w,x,y,z)$
  is totally reflexive.
\end{rmk}

\section{Families of indecomposable modules of different size} 
\label{sec:bt}

\noindent
With appropriately chosen ring elements as input, \conref{indec}
yields the infinite families of modules in the Brauer--Thrall theorems
advertised in \secref{intro}. In this section we begin to analyze the
requirements on the input.

\begin{thm}
  \label{thm:indectr}
  Let $\Rm$ be a local ring and assume that $w$ and $x$ are elements
  in $\mnotmm$, that form an exact pair of zero divisors. Assume
  further that $y$ and $z$ are elements in $\mnotmm$ with $yz=0$ and
  that one of the following conditions holds:
  \begin{prt}
  \item The elements $w$, $x$, and $y$ are linearly independent modulo
    $\m^2$.
  \item One has $w \in (x) + \m^2$ and $y,z \not \in (x)+ \m^2$.
  \end{prt}
  For every $n\in\NN$, the $R$-module $M_n(w,x,y,z)$ is
  indecomposable, totally reflexive, and non-free.  Moreover,
  $M_n(w,x,y,z)$ has constant Betti numbers, equal to~$n$, and its
  minimal free resolution is periodic of period at most $2$.
\end{thm}

The proof of \thmref[]{indectr}---which takes up the balance of the
section---employs an auxiliary result of independent interest,
\prpref{indec} below; its proof is deferred to the end of the section.

\begin{prp}
  \label{prp:indec}
  Let $\Rm$ be a local ring, let $n$ be a positive integer, and let
  $w$, $x$, $y$, and $z$ be elements in $\mnotmm$.
  \begin{prt}
  \item Assume that $w$, $x$, and $y$ are linearly independent modulo
    $\m^2$.
    \begin{itemize}
    \item If $n$ is even, then $M_n(w,x,y,z)$ is indecomposable.
    \item If $n$ is odd, then $M_n(w,x,y,z)$ or $M_n(x,w,y,z)$ is
      indecomposable.
    \end{itemize}
  \item If $y\notin (w,x) + \m^2$ and $z\notin (x) + \m^2$ hold, then
    $M_n(w,x,y,z)$ is indecomposable.
  \end{prt}
\end{prp}

\begin{prf*}[Proof of Theorem \pgref{thm:indectr}]
  Let $n$ be a positive integer; in view of \rmkref{del1}, all we need
  to show is that the $R$-module $M_n(w,x,y,z)$ is indecomposable.

  (a): Assume that $w$, $x$, and $y$ are linearly independent modulo
  $\m^2$.  By \rmkref{del1} the module $M_n(w,x,y,z)$ is the first
  syzygy of $M_n(x,w,-y,-z)$.  If the module $M_n(x,w,y,z)$ is
  indecomposable, then so is the isomorphic module $M_n(x,w,-y,-z)$,
  and it follows from \lemref{indec} that $M_n(w,x,y,z)$ is
  indecomposable as well. Thus by \partprpref{indec}{a} the module
  $M_n(w,x,y,z)$ is indecomposable.

  (b): Under the assumptions $w \in (x)+\m^2$ and $y,z\notin
  (x)+\m^2$, the conditions in \partprpref{indec}{b} are met, so the
  $R$-module $M_n(w,x,y,z)$ is indecomposable.
\end{prf*}

\begin{prf*}[Proof of Proposition \pgref{prp:indec}]
  Let $n$ be a positive integer and let $w$, $x$, $y$, and $z$ be
  elements in $\mnotmm$. It is convenient to work with a presentation
  matrix $\Phi_n(w,x,y,z)$ for $M = M_n(w,x,y,z)$ that one obtains as
  follows.  Set $p=\lceil\frac{n}{2}\rceil$ and let $\Pi$ be the
  $n\times n$ matrix obtained from $\Iota_n$ by permuting its rows
  according to
  \begin{equation*}
    \begin{pmatrix}
      1 & 2 & 3 & 4 & 5 & 6 & \cdots &n\\
      1 & p+1 & 2 & p+2 & 3 & p+3 & \dots & \delta p + p
    \end{pmatrix},
  \end{equation*}
  with $\delta=0$ if $n$ is odd and $\delta = 1$ if $n$ is even. Set
  $\Phi_n(w,x,y,z) = \Pi\Theta_n(w,x,y,z)\Pi^{-1}$.  If $n$ is even,
  then $\Phi_n$ is the block matrix
  \begin{equation*}
    \Phi_n(w,x,y,z) =
    \begin{pmatrix}
      w\Iota_p & y\Iota_p\\
      z\J_p & x\Iota_p
    \end{pmatrix},
  \end{equation*}
  and if $n$ is odd, then $\Phi_n(w,x,y,z)$ is the matrix obtained
  from $\Phi_{n+1}(w,x,y,z)$ by deleting the last row and the last
  column.

  To verify that $M$ is indecomposable, assume that $\varepsilon \in
  \Hom{M}{M}$ is idempotent and not the identity map $1_M$.  The goal
  is to show that $\varepsilon$ is the zero map.  The only idempotent
  automorphism of $M$ is $1_M$, so $\varepsilon$ is not an
  isomorphism.  Thus $\varepsilon$ is not surjective, as $M$ is
  noetherian.  Set $\Phi = \Phi_n(w,x,y,z)$ and consider the
  commutative diagram with exact rows
  \begin{equation}
    \label{eq:lift}
    \begin{gathered}
      \xymatrixcolsep{1.5pc} \xymatrix{ R^n \ar[r]^ \Phi \ar[d]^\Beta
        & R^n \ar[r] \ar[d]^\Alpha & M \ar[r] \ar[d]^\varepsilon & 0
        \\
        R^n \ar[r]^\Phi & R^n \ar[r] & M \ar[r] & 0 }
    \end{gathered}
  \end{equation}
  obtained by lifting $\varepsilon$. Let $\bar{\Alpha} = (a_{ij})$ and
  $\widebar{\Beta}$ be the $n \times n$ matrices obtained from
  $\Alpha$ and $\Beta$ by reducing their entries modulo $\m$.

  To prove that $\varepsilon$ is the zero map, it suffices to show
  that the matrix $\bar{\Alpha}$ is nilpotent. Indeed, if
  $\bar{\Alpha}$ is nilpotent, then so is the map $\bar \varepsilon
  \colon M/{\m M} \to M/{\m M}$. As $\bar{\varepsilon}$ is also
  idempotent, it is the zero map. Now it follows from Nakayama's lemma
  that the map $1_M- \varepsilon$ is surjective and hence an
  isomorphism. As $1_M- \varepsilon$ is idempotent, it follows that
  $1_M - \varepsilon$ is the identity map $1_M$; that is, $\varepsilon
  = 0$.

  \subparagraph \emph{Claim.} If the matrix $\bar{\Alpha}$ is upper
  triangular with identical entries on the diagonal, i.e.~$a_{11} =
  a_{22} = \cdots = a_{nn}$, then it is nilpotent.

  \emph{Proof.} Since $\varepsilon$ is not surjective, the matrix
  $\Alpha$ does not represent a surjective map and, by Nakayama's
  lemma, neither does $\bar{\Alpha}$. Therefore, the diagonal entries
  of $\bar{\Alpha}$ cannot all be non-zero, whence they are all zero,
  and $\bar{\Alpha}$ is nilpotent.

  \subparagraph Denote by $\w$, $\x$, $\y$, and $\z$ the images of
  $w$, $x$, $y$, and $z$ in $\m/\m^2$, and let $V$ be the
  $\k$-subspace of $\m/\m^2$ spanned by $\w$, $\x$, $\y$, and $\z$.
  Consider the following possibilities:
  \begin{itemize}
  \item[(I)] The elements $\w$, $\x$, $\y$, and $\z$ form a basis for
    $V$.
  \item[(II)] The elements $\x$, $\y$, and $\z$ form a basis for $V$,
    and $\k\w= \k\x$ holds.
  \item[(III)] The elements $\w$, $\x$, and $\y$ form a basis for $V$.
  \item[(IV)] The elements $\x$ and $\y$ form a basis for $V$, and one
    has $\k\w=\k\x$ and $\z\notin \k\x$.
  \end{itemize}
  Under the assumptions on $w$, $x$, and $y$ in part (a), one of the
  conditions (I) or (III) holds. Under the assumptions in part (b),
  one of the conditions (I)--(IV) holds. Indeed, if $V$ has dimension
  $4$, then (I) holds. If that is not the case, then the dimension of
  $V$ is $2$ or $3$.  In case $\dim[\k]{V} = 2$, the elements $\x$ and
  $\y$ form a basis for $V$, whereas $\w$ and $\x$ cannot be linearly
  independent; thus (IV) holds.  In case $\dim[\k]{V} = 3$, condition
  (II) or (III) holds, depending on whether or not the equality $\k\w
  = \k\x$ holds.

  The rest of the proof is split in two, according to the parity of
  $n$. To prove that $M_n(w,x,y,z)$ is indecomposable, it suffices to
  prove that the matrix $\bar{\Alpha}$ is nilpotent. This is how we
  proceed under each of the conditions (I), (II), and (IV), and under
  condition (III) when $n$ is even. When $n$ is odd, and condition
  (III) holds, we show that one of the modules $M_n(w,x,y,z)$ and
  $M_n(x,w,y,z)$ is indecomposable.

  \subparagraph \emph{Case 1: $n$ is even.}  Let $\widetilde{\Phi}$
  denote the matrix obtained from $\Phi$ by reducing the entries
  modulo $\m^2$. Write $\bar{\Alpha}$ and $\widebar{\Beta}$ as block
  matrices
  \begin{equation*}
    \bar{\Alpha}=
    \begin{pmatrix}
      \Alpha_{11}& \Alpha_{12}\\
      \Alpha_{21} & \Alpha_{22}
    \end{pmatrix}
    \qqtext{and} \widebar{\Beta}=
    \begin{pmatrix}
      \Beta_{11}& \Beta_{12}\\
      \Beta_{21} & \Beta_{22}
    \end{pmatrix},
  \end{equation*}
  where $\Alpha_{ij}$ and $\Beta_{ij}$ are $p\times p$ matrices with
  entries from $\k$.  By \eqref{lift}, the equality $\Alpha\Phi =
  \Phi\Beta$ holds; it implies an equality of block matrices
  \begin{equation}
    \label{eq:relations}
    \begin{pmatrix}
      \w \Alpha_{11} + \z \Alpha_{12}\J_p & \y \Alpha_{11} + \x \Alpha_{12}\\
      \w \Alpha_{21}+\z \Alpha_{22}\J_p & \y \Alpha_{21}+\x
      \Alpha_{22}
    \end{pmatrix}
    =
    \begin{pmatrix}
      \w \Beta_{11} + \y \Beta_{21} & \w \Beta_{12} + \y \Beta_{22}\\
      \z \J_p\Beta_{11} + \x \Beta_{21} & \z \J_p\Beta_{12} + \x
      \Beta_{22}
    \end{pmatrix}.
  \end{equation}

  Assume first that condition (I) or (III) holds, so that the elements
  $\w$, $\x$, and $\y$ are linearly independent. Then the equality of
  the upper right blocks, $\y \Alpha_{11} + \x \Alpha_{12} = \w
  \Beta_{12} + \y \Beta_{22}$, yields
  \begin{equation*}
    \Alpha_{11} = \Beta_{22} \qqand \Alpha_{12} = 0 = \Beta_{12}.
  \end{equation*}
  From the blocks on the diagonals, one now gets
  \begin{equation*}
    \Alpha_{11} = \Beta_{11}, \qquad \Alpha_{21} = 0 = \Beta_{21}, \qqand \Alpha_{22} =
    \Beta_{22}.
  \end{equation*}
  Thus the matrix $\bar{\Alpha}$ has the form
  $\left(\begin{smallmatrix}\Alpha_{11} & 0\\ 0 &
      \Alpha_{11}\end{smallmatrix}\right)$.  Finally, the equality of
  the lower left blocks yields $\Alpha_{11}\J_p = \J_p\Alpha_{11}$.
  Since $\J_p$ is non-derogatory, this implies that $\Alpha_{11}$
  belongs to the algebra $\k[\J_p]$ of polynomials in $\J_p$. That is,
  there are elements $c_0,\ldots,c_{p-1}$ in $\k$ such that
  $\Alpha_{11}=c_0\Iota_p + c_1\J_p + \cdots + c_{p-1}\J_p^{p-1}$; see
  \thmcite[3.2.4.2]{matrixanalysis}.  In particular, the matrices
  $\Alpha_{11}$ and, therefore, $\bar{\Alpha}$ are upper-triangular
  with identical entries on the diagonal. By Claim, $\bar{\Alpha}$ is
  nilpotent as desired.

  Under either condition (II) or (IV), the elements $\x$ and $\y$ are
  linearly independent, $\z$ is not in $\k\x$, and the equality $\k\w
  = \k\x$ holds.  In particular, there is an element $t \ne 0$ in $\k$
  such that $t\w = \x$.  From the off-diagonal blocks in
  \eqref{relations} one obtains the following relations:
  \begin{equation}
    \label{eq:gurk}
    \begin{split}
      \begin{aligned}
        \Alpha_{11} &= \Beta_{22}\\ t\Alpha_{12} &= \Beta_{12}
      \end{aligned}
      \qqand
      \begin{aligned}
        \Alpha_{22}\J_p &= \J_p\Beta_{11} \\ \Alpha_{21} &=
        t\Beta_{21}.
      \end{aligned}
    \end{split}
  \end{equation}

  If (II) holds, then the blocks on the diagonals in \eqref{relations}
  yield
  \begin{equation*}
    \Alpha_{11} = \Beta_{11},\qquad \Alpha_{22} = \Beta_{22},\qqand \Alpha_{21} = 0.
  \end{equation*}
  Thus the matrix $\bar{\Alpha}$ has the form
  $\left(\begin{smallmatrix}\Alpha_{11} & \Alpha_{12}\\ 0 &
      \Alpha_{11}\end{smallmatrix}\right)$, and the equality
  $\Alpha_{11}\J_p = \J_p\Alpha_{11}$ holds. As above, it follows that
  the matrices $\Alpha_{11}$ and, therefore, $\bar{\Alpha}$ are
  upper-triangular with identical entries on the diagonal. That is,
  $\bar{\Alpha}$ is nilpotent by Claim.

  Now assume that (IV) holds.  There exist elements $r$ and $s \ne 0$
  in $\k$ such that $\z = r\x + s\y$. Comparison of the blocks on the
  diagonals in \eqref{relations} now yields
  \begin{equation}
    \label{eq:morequrk}
    \begin{split}
      \begin{aligned}
        \Alpha_{11} + rt\Alpha_{12}\J_p &= \Beta_{11}\\
        s \Alpha_{12}\J_p &= \Beta_{21}
      \end{aligned}
      \qqand
      \begin{aligned}
        \Alpha_{22} &=r \J_p\Beta_{12} + \Beta_{22}\\
        \Alpha_{21} &= s \J_p \Beta_{12}.
      \end{aligned}
    \end{split}
  \end{equation}
  Combine these equalities with those from \eqref{gurk} to get
  \begin{equation*}
    \bar{\Alpha} =
    \begin{pmatrix}
      \Alpha_{11} & \Alpha_{12}\\
      st \J_p\Alpha_{12} & \Alpha_{11} + rt\J_p\Alpha_{12}
    \end{pmatrix}.
  \end{equation*}
  It follows from the equalities
  \begin{equation*}
    \J_p\Alpha_{12} = t^{-1}\J_p\Beta_{12} = (st)^{-1}\Alpha_{21}
    =s^{-1}\Beta_{21} = \Alpha_{12}\J_p,
  \end{equation*}
  derived from \eqref{gurk} and \eqref{morequrk}, that the matrix
  $\Alpha_{12}$ commutes with $\J_p$; hence it belongs to $\k[\J_p]$.
  Similarly, the chain of equalities
  \begin{align*}
    \J_p\Alpha_{11} = \J_p\Beta_{11} - rt\J_p\Alpha_{12}\J_p =&
    \J_p\Beta_{11} - r \J_p\Beta_{12}\J_p\\ =& \Alpha_{22}\J_p - r
    \J_p\Beta_{12}\J_p = \Beta_{22}\J_p = \Alpha_{11}\J_p
  \end{align*}
  shows that $\Alpha_{11}$ is in $\k[\J_p]$. Thus all four blocks in
  $\bar{\Alpha}$ belong to $\k[\J_p]$.  For notational bliss, identify
  $\k[\J_p]$ with the ring $S = \k[\chi]/(\chi^p)$, where $\chi$
  corresponds to $t\J_p$. With this identification, $\bar{\Alpha}$
  takes the form of a $2 \times 2$ matrix with entries in $S$:
  \begin{equation*}
    \begin{pmatrix} f & g \\ sg\chi & f + rg\chi \end{pmatrix}.
  \end{equation*}
  As $\bar{\Alpha}$ is not invertible, the determinant $f^2 + frg\chi
  - sg^2\chi$ belongs to the maximal ideal $(\chi)$ of $S$.  It
  follows that $f$ is in $(\chi)$, whence one has
  $\bar{\Alpha}^{2p}=0$ as desired.

  \subparagraph \emph{Case 2: n is odd.}  Set $q=p-1$, where $p =
  \lceil\frac{n}{2}\rceil = \frac{n+1}{2}$.  The presentation matrix
  $\Phi$ takes the form
  \begin{equation*}
    \Phi =
    \begin{pmatrix}
      w\Iota_p & y\Eta \\
      z\Kappa & x\Iota_q
    \end{pmatrix},
  \end{equation*}
  where $\Eta$ and $\Kappa$ are the following block matrices $\Eta =
  \smash{\begin{pmatrix} \Iota_q\\0_{1\times q}\end{pmatrix}}$ and
  $\Kappa = \begin{pmatrix} 0_{q\times
      1}&\mspace{-8mu}\Iota_q\end{pmatrix}$. Notice that there are
  equalities
  \begin{equation}
    \label{eq:zero}
    \Eta\Chi = \begin{pmatrix} \Chi \\ 0_{1\times m} \end{pmatrix}
    \qqand
    \Chi'\Kappa = \begin{pmatrix} 0_{m'\times 1} &\mspace{-8mu}\Chi'\end{pmatrix}
  \end{equation}
  for every $q\times m$ matrix $\Chi$ and every $m'\times q$ matrix
  $\Chi'$.  Furthermore, it is straightforward to verify the
  equalities
  \begin{equation}
    \label{eq:hk}
    \Eta\Kappa = \J_p \qqand \Kappa\Eta = \J_q.
  \end{equation}
  As in Case 1, write
  \begin{equation*}
    \bar{\Alpha}=
    \begin{pmatrix}
      \Alpha_{11}& \Alpha_{12}\\
      \Alpha_{21} & \Alpha_{22}
    \end{pmatrix}
    \qqtext{and}
    \widebar{\Beta}=
    \begin{pmatrix}
      \Beta_{11}& \Beta_{12}\\
      \Beta_{21} & \Beta_{22}
    \end{pmatrix},
  \end{equation*}
  where, now, $\Alpha_{ij}$ and $\Beta_{ij}$ are matrices of size
  $m_i\times m_j$, for $m_1=p$ and $m_2=q$.  With $\widetilde{\Phi}$
  as defined in Case 1, the relation $\bar{\Alpha}\widetilde{\Phi} =
  \widetilde{\Phi}\widebar{\Beta}$, derived from \eqref{lift}, yields:
  \begin{equation}
    \label{eq:oddrelations}
    \begin{pmatrix}
      \w \Alpha_{11} + \z \Alpha_{12}\Kappa & \y \Alpha_{11}\Eta + \x \Alpha_{12}\\
      \w \Alpha_{21} + \z \Alpha_{22}\Kappa & \y \Alpha_{21}\Eta + \x
      \Alpha_{22}
    \end{pmatrix}
    =
    \begin{pmatrix}
      \w \Beta_{11} + \y \Eta\Beta_{21} & \w \Beta_{12} + \y \Eta\Beta_{22}\\
      \z \Kappa\Beta_{11} + \x \Beta_{21} & \z \Kappa\Beta_{12} + \x
      \Beta_{22}
    \end{pmatrix}.%
  \end{equation}

  Assume first that condition (I) or (III) holds, so that the elements
  $\w$, $\x$, and $\y$ are linearly independent. From the equality of
  the upper right blocks in \eqref{oddrelations} one~gets
  \begin{equation}
    \label{eq:purk}
    \Alpha_{11}\Eta = \Eta\Beta_{22} \qqand \Alpha_{12} = 0 = \Beta_{12}.
  \end{equation}
  In view of \eqref{zero}, comparison of the blocks on the diagonals
  now yields
  \begin{equation}
    \label{eq:odd1}
    \begin{split}
      \begin{aligned}
        \Alpha_{11} &= \Beta_{11}\\ \Beta_{21} &= 0
      \end{aligned}
      \qqand
      \begin{aligned}
        \Alpha_{21}\Eta &= 0\\ \Alpha_{22} &= \Beta_{22}.
      \end{aligned}
    \end{split}
  \end{equation}
  From the first equality in \eqref{purk} and the last equality in
  \eqref{odd1} one gets in view of \eqref{zero}
  \begin{equation}
    \label{eq:alpha-info}
    \Alpha_{11} =
    \begin{pmatrix}
      \Alpha_{22} & *\\
      0_{1\times q} & *
    \end{pmatrix};
  \end{equation}
  here and in the following the symbol `$*$' in a matrix denotes an
  unspecified block of appropriate size.  To glean information from
  the equality of the lower left blocks in \eqref{oddrelations},
  assume first that $\w$ and $\z$ are linearly independent.  Then one
  has
  \begin{equation}
    \label{eq:odd2}
    \Alpha_{21} = 0 \qqand \Alpha_{22}\Kappa = \Kappa\Beta_{11}.
  \end{equation}
  Combine this with \eqref{alpha-info} and the second equality in
  \eqref{purk} to see that the matrix $\bar{\Alpha}$ has the form
  \begin{equation}
    \label{eq:aform}
    \bar{\Alpha} =
    \left(
      \begin{array}{c|c}
        \begin{matrix} \Alpha_{22} & *\\ 0_{1\times q} & * \end{matrix}
        & 0_{p\times q}\\ \hline
        0_{q\times p}  & \Alpha_{22}
      \end{array}
    \right).
  \end{equation}
  The equalities
  \begin{equation*}
    \Alpha_{11}\J_p = \Alpha_{11}\Eta\Kappa = \Eta\Beta_{22}\Kappa =
    \Eta\Alpha_{22}\Kappa =  \Eta\Kappa\Beta_{11} = \J_p\Alpha_{11}
  \end{equation*}
  and
  \begin{equation*}
    \Alpha_{22}\J_q = \Alpha_{22}\Kappa\Eta = \Kappa\Beta_{11}\Eta =
    \Kappa\Alpha_{11}\Eta =  \Kappa\Eta\Beta_{22} = \J_q\Alpha_{22} ,
  \end{equation*}
  derived from \eqref{hk}, \eqref{purk}, \eqref{odd1}, and
  \eqref{odd2}, show that $\Alpha_{11}$ is in $\k[\J_p]$ and
  $\Alpha_{22}$ is in $\k[\J_q]$. It follows that $\bar{\Alpha}$ is
  upper triangular with identical entries on the diagonal. Thus
  $\bar{\Alpha}$ is nilpotent, and $M_n(w,x,y,z)$ is
  indecomposable. If, on the other hand, $\z$ and $\w$ are linearly
  dependent, then $\z$ and $\x$ are linearly independent, as $\w$ and
  $\x$ are linearly independent by assumption.  It follows from what
  we have just shown that $M_n(x,w,y,z)$ is indecomposable.

  Under either condition (II) or (IV), the elements $\x$ and $\y$ are
  linearly independent, $\z$ is not in $\k\x$, and the equality $\k\w
  = \k\x$ holds.  In particular, there is an element $t \ne 0$ in $\k$
  such that $t\w = \x$. Compare the off-diagonal blocks in
  \eqref{oddrelations} to get
  \begin{equation}
    \label{eq:odd3}
    \begin{split}
      \begin{aligned}
        t\Alpha_{12} &= \Beta_{12}\\ \Alpha_{11}\Eta &= \Eta\Beta_{22}
      \end{aligned}
      \qqand
      \begin{aligned}
        \Alpha_{21} &= t\Beta_{21}\\ \Alpha_{22}\Kappa &=
        \Kappa\Beta_{11}.
      \end{aligned}
    \end{split}
  \end{equation}

  If (II) holds, then a comparison of the blocks on the diagonals in
  \eqref{oddrelations} combined with \eqref{zero} implies
  \begin{equation}
    \label{eq:CaseIIodd}
    \Alpha_{11} = \Beta_{11}, \quad \Alpha_{12} = 0 = \Beta_{21}, \qand \Alpha_{22} =
    \Beta_{22}.
  \end{equation}
  It follows from \eqref{odd3} and \eqref{CaseIIodd} that the matrix
  $\Alpha_{21}$ is zero. In view of the equality $\Alpha_{11}\Eta =
  \Eta\Beta_{22}$ from \eqref{odd3}, it now follows that
  $\bar{\Alpha}$ has the form given in \eqref{aform}.  Using the
  equalities in \eqref{odd3} and \eqref{CaseIIodd}, one can repeat the
  arguments above to see that $\Alpha_{11}$ is in $\k[\J_p]$ and
  $\Alpha_{22}$ is in $\k[\J_q]$, and continue to conclude that
  $\bar{\Alpha}$ is nilpotent.

  Finally, assume that (IV) holds.  There exist elements $r$ and $s
  \ne 0$ in $\k$ such that $\z = r\x + s\y$. Comparison of the blocks
  on the diagonals in \eqref{oddrelations} yields
  \begin{equation}
    \label{eq:oddgurk}
    \begin{split}
      \begin{aligned}
        \Alpha_{11} + rt\Alpha_{12}\Kappa &= \Beta_{11}\\
        s \Alpha_{12}\Kappa &= \Eta\Beta_{21}
      \end{aligned}
      \qqand
      \begin{aligned}
        \Alpha_{21}\Eta &=s \Kappa\Beta_{12}\\
        \Alpha_{22} &= \Beta_{22} + r \Kappa\Beta_{12}.
      \end{aligned}
    \end{split}
  \end{equation}
  The equality $s \Alpha_{12}\Kappa = t^{-1} \Eta\Alpha_{21}$,
  obtained from \eqref{odd3} and \eqref{oddgurk}, shows, in view of
  \eqref{zero}, that the matrices $\Alpha_{12}$ and $\Alpha_{21}$ have
  the following form:
  \begin{equation*}
    \Alpha_{12} = \begin{pmatrix} * \\ \Epsilon \end{pmatrix}
    \qqand
    \Alpha_{21} = \begin{pmatrix} 0_{q\times 1} &\mspace{-12mu} \Gamma \end{pmatrix},
  \end{equation*}
  where $\Epsilon$ and $\Gamma$ are $q \times q$ matrices, and the
  last row of $\Epsilon$ is zero.  From the equalities in \eqref{hk},
  \eqref{odd3}, and \eqref{oddgurk} one gets
  \begin{equation*}
    \Alpha_{21}\J_p= s \Kappa\Beta_{12}\Kappa = st\Kappa\Alpha_{12}\Kappa =
    t\Kappa\Eta\Beta_{21} = \J_q\Alpha_{21}.
  \end{equation*}
  From here it is straightforward to verify that $\Gamma$ commutes
  with $\J_q$; i.e.\ $\Gamma$ belongs to $\k[\J_q]$. Similarly, from
  the equalities
  \begin{equation*}
    \Alpha_{12}\J_q = s^{-1}\Eta\Beta_{21}\Eta = (st)^{-1} \Eta\Alpha_{21}\Eta =
    t^{-1} \Eta\Kappa\Beta_{12} = \J_p\Alpha_{12},
  \end{equation*}
  it follows that $\Epsilon$ belongs to $\k[\J_q]$. Since the last row
  in $\Epsilon$ is zero, all entries on the diagonal of $\Epsilon$ are
  zero, and the matrix is nilpotent. The first equality in
  \eqref{oddgurk} can now be written as
  \begin{equation*}
    \Alpha_{11} = \Beta_{11} - \begin{pmatrix} 0 & * \\ 0_{q\times 1} & rt\Epsilon \end{pmatrix}.
  \end{equation*}
  Combine this with the last equality in \eqref{odd3} to get
  \begin{equation*}
    \bar{\Alpha} =
    \left(
      \begin{array}{cc|c}
        a_{11} & * & *\\
        0_{q\times 1} & \Alpha_{22} - rt\Epsilon & \Epsilon\\ \hline
        0_{q\times 1}& \Gamma  & \Alpha_{22}
      \end{array}
    \right).
  \end{equation*}
  The equalities
  \begin{align*}
    \Alpha_{11}\J_p = \Eta\Beta_{22}\Kappa =\Eta(\Alpha_{22} - r
    \Kappa\Beta_{12})\Kappa &= \Eta\Alpha_{22}\Kappa - r
    \Eta\Kappa\Beta_{12}\Kappa\\ &= \Eta\Kappa\Beta_{11} -
    rt\Eta\Kappa\Alpha_{12}\Kappa = \J_p\Alpha_{11}
  \end{align*}
  show that $\Alpha_{11}$ is in $\k[\J_p]$, and a similar chain of
  equalities shows that $\Alpha_{22}$ is in $\k[\J_q]$. It follows
  that all entries on the diagonal of $\bar{\Alpha}$ are
  identical. Let $\Delta$ be the matrix obtained by deleting the first
  row and first column in $\bar{\Alpha}$ and write it in block form
  \begin{equation*}
    \Delta=
    \begin{pmatrix}
      \Alpha_{22} - rt \Epsilon & \Epsilon\\
      \Gamma & \Alpha_{22}
    \end{pmatrix}.
  \end{equation*}
  As $\bar{\Alpha}$ is not invertible, one has $0=\det{\bar{\Alpha}} =
  a_{11}(\det{\Delta})$. If $a_{11}$ is non-zero, then $\Delta$ has
  determinant $0$; in particular, it is not invertible. Considered as
  a $2\times 2$ matrix over the artinian local ring $\k[\J_q]$, its
  determinant $(\Alpha_{22})^2 - rt\Epsilon\Alpha_{22} -
  \Gamma\Epsilon$ belongs to the maximal ideal $(\J_q)$. As $\Epsilon$
  is nilpotent, it belongs to $(\J_q)$ and hence so does
  $\Alpha_{22}$; this contradicts the assumption that $a_{11}$ is
  non-zero. If the diagonal entry $a_{11}$ is $0$, then the matrix
  $\Alpha_{22}$ is nilpotent, which implies that $\Delta$ is nilpotent
  and, finally, that $\bar{\Alpha}$ is nilpotent.
\end{prf*}

\section{Brauer--Thrall I over short local rings with exact zero
  divisors} 
\label{sec:m3}

\noindent
Let $\Rmk$ be a local ring. The \emph{embedding dimension} of $R$,
denoted $\embdim R$, is the minimal number of generators of $\m$,
i.e.\ the dimension of the $\k$-vector space $\m/\m^2$. The
\emph{Hilbert series} of $R$ is the power series
$\operatorname{H}_R(\tau) = \sum_{i=0}^\infty
\dim[\k]{(\m^{i}/\m^{i+1})}\tau^i$.

In the rest of this section $\Rm$ is a local ring with $\m^3=0$. The
main result, \thmref{m3zd}, together with \eqref{length} and
\thmref{indectr}, establishes \thmref{firstbt}. Towards the proof of
\thmref{m3zd}, we first recapitulate a few facts about totally
reflexive modules and exact zero divisors.

If $R$ is Gorenstein, then every $R$-module is totally reflexive; see
\thmcite[(4.13) and (4.20)]{MAsMBr69}. If $R$ is not Gorenstein, then
existence of a non-free totally reflexive $R$-module forces certain
relations among invariants of $R$. The facts in \pgref{m3tr} are
proved by Yoshino \cite{YYs03}; see also \cite{LWCOVl07} for the
non-graded case.

\begin{bfhpg}[Totally reflexive modules]
  \label{m3tr}
  Assume that $R$ is not Gorenstein and set $e=\embdim R$.  If $M$ is
  a totally reflexive $R$-module without free summands and minimally
  generated by $n$ elements, then the equalities
  \begin{equation}
    \label{eq:length}
    \lgt{M} =ne  \qqand \bet{i}{M} = n \ \text{ for all $i\ge 0$}
  \end{equation}
  hold. Moreover, $\m^2$ is non-zero and the following hold:
  \begin{equation}
    \label{eq:m3tr}
    \ann{\m} = \m^2, \quad  \dim[\k]{\m^2} = e-1, \qand  \lgtR = 2e.
  \end{equation}
  In particular, $e$ is at least $3$, and the Hilbert series of $R$ is
  $1 + e\tau +(e-1)\tau^2$.
\end{bfhpg}

Let $\k$ be a field. For the Gorenstein ring $R=\k[x]/(x^3)$, two of
the relations in \eqref{m3tr} fail. This ring also has an exact zero
divisor, $x^2$, in the square of the maximal ideal; for rings with
embedding dimension $2$ or higher this cannot happen.

\begin{bfhpg}[Exact zero divisors]
  \label{m3epzd}
  Set $e=\embdim R$ and assume $e\ge 2$.  Suppose that $w$ and $x$
  form an exact pair of zero divisors in $R$. As $\m^2$ is contained
  in the annihilator of $\m$, there is an inclusion $\m^2\subseteq
  (x)$, which has to be strict, as $(w) = \ann{x}$ is strictly
  contained in $\m$. Thus $x$ is a minimal generator of $\m$ with
  \begin{equation*}
    x\m = \m^2.
  \end{equation*}
  By symmetry, $w$ is a minimal generator of $\m$ with $w\m = \m^2$.
  Let $\set{v_1,\ldots,v_{e-1},w}$ be a minimal set of generators for
  $\m$, then the elements $xv_1,\dots,xv_{e-1}$ generate $\m^2$, and
  it is elementary to verify that they form a basis for $\m^2$ as a
  $\k$-vector space.  It follows that the relations in \eqref{m3tr}
  hold and, in addition, there is an equality
  \begin{equation*}
    \lgt{(x)} = e.
  \end{equation*}
  Note that the socle $\ann{\m}$ of $R$ has dimension $e-1$ over $\k$,
  so $R$ is Gorenstein if and only if $e=2$.
\end{bfhpg}

\begin{lem}
  \label{lem:ann}
  Assume that $\Rm$ has Hilbert series $1 + e\tau +(e-1)\tau^2$ with
  $e \ge 2$. For every element $x \in \mnotmm$ the following hold:
  \begin{prt}
  \item The ideal $(x)$ in $R$ has length at most $e$.
  \item There exists an element $w\in\mnotmm$ that annihilates $x$,
    and if $w$ generates $\ann{x}$, then $w$ and $x$ form an exact
    pair of zero divisors in~$R$.
  \item If the equalities $wx=0$ and $w\m = \m^2 = x\m$ hold, then $w$
    and $x$ form an exact pair of zero divisors in $R$.
  \end{prt}
\end{lem}

\begin{prf*}
  (a): By assumption, the length of $\m$ is $2e-1$.  As $x$ and $e-1$
  other elements form a minimal set of generators for $\m$, the
  inequality $2e-1\ge \lgt{(x)}+e-1$ holds; whence $\lgt{(x)}$ is at
  most $e$.

  (b): Additivity of length on short exact sequences yields
  \begin{equation*}
    \lgt{\ann{x}} = \lgtR - \lgt{(x)} \ge e,
  \end{equation*}
  so $\m^2$ is properly contained in $\ann{x}$; choose an element $w$
  in $\ann{x}\backslash \m^2$. There is an inclusion $(x)\subseteq
  \ann{w}$, and if the equality $(w) = \ann{x}$ holds, then two length
  counts yield
  \begin{equation*}
    \lgt{\ann{w}} = \lgtR-\lgt{(w)}=\lgt{(x)}.
  \end{equation*}
  Thus $(x)=\ann{w}$ holds; hence $w$ and $x$ form an exact pair of
  zero~divisors in $R$.

  (c): As both ideals $(w)$ and $(x)$ strictly contain $\m^2$, the
  equalities $\lgt{(x)} = e = \lgt{(w)}$ hold in view of part (a). As
  $w$ and $x$ annihilate each other, simple length counts show that
  they form an exact pair of zero divisors in $R$.
\end{prf*}

\begin{thm}
  \label{thm:m3zd}
  Let $\Rm$ be a local ring with $\m^3=0$ and $\embdim R \ge 3$.
  Assume that $w$ and $x$ form an exact pair of zero divisors in
  $R$. For every element $y$ in $\m\backslash(w,x)$ there exists an
  element $z\in \m\backslash\m^2$ such that the $R$-modules
  $M_n(w,x,y,z)$ are indecomposable and totally reflexive for all
  $n\in\NN$.
\end{thm}

\begin{prf*}
  By \pgref{m3epzd} the equalities in \eqref{m3tr} hold for $R$.  Let
  $y$ be an element in $\m\backslash (w,x)$. By
  \partlemref{ann}{b} the ideal $\ann{y}$ contains an element $z$ of
  $\m\backslash \m^2$.  Since $y$ is not contained in the ideal $(w,x)
  = (w,x) + \m^2$, the element $z$ is not in $(x) = (x) + \m^2$. The
  desired conclusion now follows from \thmref{indectr}.
\end{prf*}

The key to the theorem is that existence of an exact zero divisor in
$R$ implies the existence of additional elements such that the
conditions in \thmref{indectr} are satisfied. This phenomenon does not
extend to rings with $\m^4=0$.

\begin{exa}
  Let $\mathsf{F}$ be a field and set
  $S=\poly[\mathsf{F}]{x,y,z}/(x^2,y^2z,yz^2,y^3,z^3)$; it is a
  standard graded $\mathsf{F}$-algebra with Hilbert series $1 + 3\tau
  + 5\tau^2 + 3\tau^3$, and $x$ is an exact zero divisor in $S$. Set
  $\n=(x,y,z)S$ and let $v$ be an element in $\n\backslash
  ((x)+\n^2)$. A straightforward calculation shows that the
  annihilator $\ann{v}$ is contained in $\n^2$.
\end{exa}

\section{Exact zero divisors from totally reflexive modules} 
\label{sec:exist}

\noindent
Let $\Rm$ be a local ring with $\m^3=0$. If $R$ is not Gorenstein,
then a cyclic totally reflexive $R$-module is either free or generated
by an exact zero divisor. Indeed, if it is not free, then by
\eqref{length} it has constant Betti numbers, equal to $1$, so by
\rmkref{ezdtac} it is generated by an exact zero divisor in $R$. The
next results improve on this elementary observation; in particular,
\corref{diagonal} should be compared to \rmkref{ezdtac}.

\begin{lem}
  \label{lem:pm}
  Let $\Rm$ be a local ring with $\m^3=0$ and let
  \begin{equation*}
    F : \quad  F_2 \lra F_1 \xra{\psi} F_0 \xra{\varphi} F_{-1}
  \end{equation*}
  be an exact sequence of finitely generated free $R$-modules, where
  the homomorphisms are represented by matrices with entries in
  $\m$. Let $\Psi$ be any matrix that represents $\psi$. For every row
  $\Psi_r$ of $\Psi$ the following hold:
  \begin{prt}
  \item The ideal $\mathfrak{r}$, generated by the entries of
    $\Psi_r$, contains $\m^2$.
  \item If $\,\dim[\k]{\m^2}$ is at least $2$ and $\Hom{F}{R}$ is
    exact, then $\Psi_r$ has an entry from $\mnotmm$, the entries in
    $\Psi_r$ from $\mnotmm$ generate $\mathfrak{r}$, and
    $\m\mathfrak{r} = \m^2$ holds.
  \end{prt}
\end{lem}

\begin{prf*}
  Let $\Psi$ and $\Phi$ be the matrices for $\psi$ and $\phi$ with
  respect to bases $\mathcal{B}_1$, $\mathcal{B}_0$, and
  $\mathcal{B}_{-1}$ for $F_1$, $F_0$, and $F_{-1}$. For every $p \ge
  1$, let $\set{e_1,\ldots,e_p}$ be the standard basis for $R^p$.  The
  matrix $\Phi$ is of size $l \times m$, and $\Psi$ is of size
  $m\times n$, where $l$, $m$ and $n$ denote the ranks of $F_{-1}$,
  $F_0$, and $F_1$, respectively.  We make the identifications $F_{-1}
  = R^l$, $F_0 = R^m$, and $F_1 = R^n$, by letting $\mathcal{B}_{-1}$,
  $\mathcal{B}_0$, and $\mathcal{B}_1$ correspond to the standard
  bases.  The map $\psi$ is now left multiplication by the matrix
  $\Psi$, and $\phi$ is left multiplication by $\Phi$.  For every
  $x\in \m^2$ and every basis element $e_i$ in $R^m$ one has
  $\varphi(xe_i) = 0$; indeed, $\Phi$ has entries in $\m$, the entries
  of $xe_i$ are in $\m^2$, and $\m^3=0$ holds by assumption.  By
  exactness of $F$, the element $xe_i$ is in the image of $\psi$, and
  (a) follows.

  (b): Fix $q\in\set{1,\ldots,m}$ and let $\Psi_q$ be the $q$th row of
  $\Psi = (x_{ij})$; we start by proving the following:

  \subparagraph \emph{Claim.} Every entry from $\m^2$ in $\Psi_q$ is
  contained in the ideal generated by the other entries in $\Psi_q$.

  \emph{Proof.} Assume, towards a contradiction, that some entry from
  $\m^2$ in $\Psi_q$ is not in the ideal generated by the other
  entries. After a permutation of the columns of $\Psi$, one can
  assume that the entry $x_{q1}$ is in $\m^2$ but not in the ideal
  $(x_{q2},\ldots,x_{qn})$. Since the element $x_{q1}e_q$ belongs to
  $\Ker{\varphi} = \Im{\psi}$, there exist elements $a_i$ in $R$ such
  that $\psi(\sum_{i=1}^na_ie_i) = x_{q1}e_q$.  In particular, one has
  $\sum_{i=1}^na_ix_{qi} = x_{q1}$, whence it follows that $a_1$ is
  invertible. The $n\times n$ matrix
  \begin{equation*}
    \Alpha =
    \begin{pmatrix}
      a_1 & 0 & \cdots & 0\\
      a_2 & 1 & & 0\\
      \vdots & &\ddots &\vdots\\
      a_n & 0 &\cdots & 1
    \end{pmatrix}
  \end{equation*}
  is invertible, as it has determinant $a_1$. The first column of the
  matrix $\Psi\Alpha$, which is the first row of the transposed matrix
  $(\Psi\Alpha)^\mathrm{T}$, has only one non-zero entry, namely
  $x_{q1}$.  As $\Psi\Alpha$ represents $\psi$, the matrix
  $(\Psi\Alpha)^\mathrm{T}$ represents the dual homomorphism
  $\Hom{\psi}{R}$.  By assumption the sequence $\Hom{F}{R}$ is exact,
  so it follows from part (a) that the element $x_{q1}$ spans
  $\m^2$. This contradicts the assumption that $\m^2$ is a $\k$-vector
  space of dimension at least $2$. This finishes the proof of Claim.

  \subparagraph Suppose, for the moment, that every entry of $\Psi_r$
  is in $\m^2$.  Performing column operations on $\Psi$ results in a
  matrix that also represents $\psi$, so by Claim one can assume that
  $\Psi_r$ is the zero row, which contradicts part (a).  Thus $\Psi_r$
  has an entry from $\mnotmm$. After a permutation of the columns of
  $\Psi$, one may assume that the entries $x_{r1},\ldots,x_{rt}$ are
  in $\mnotmm$ while $x_{r(t+1)}, \ldots, x_{rn}$ are in $\m^2$, where
  $t$ is in $\set{1,\ldots,n}$. Claim shows that---after column
  operations that do not alter the first $t$ columns---one can assume
  that the entries $x_{r(t+1)}, \ldots, x_{rn}$ are zero. Thus the
  entries $x_{r1},\ldots,x_{rt}$ from $\mnotmm$ generate the ideal
  $\mathfrak{r}$.

  Finally, after another permutation of the columns of $\Psi$, one can
  assume that $\set{x_{r1}, \ldots, x_{rs}}$ is maximal among the
  subsets of $\set{x_{r1}, \ldots, x_{rt}}$ with respect to the
  property that its elements are linearly independent modulo
  $\m^2$. Now use column operations to ensure that the elements
  $x_{r(s+1)},\ldots,x_{rn}$ are in $\m^2$. As above, it follows that
  $x_{r1}, \ldots, x_{rs}$ generate $\mathfrak{r}$. To verify the last
  equality in (b), note first the obvious inclusion $\m\mathfrak{r}
  \subseteq \m^2$.  For the reverse inclusion, let $x\in \m^2$ and
  write $x=x_{r1}b_1+\dots+x_{rs}b_s$ with $b_i\in R$.  If some $b_i$
  were a unit, then the linear independence of the elements $x_{ri}$
  modulo $\m^2$ would be contradicted.  Thus each $b_i$ is in $\m$,
  and the proof is complete.
\end{prf*}

The condition $\dim[\k]{\m^2} \ge 2$ in part (b) of the lemma cannot
be relaxed:

\begin{exa}
  Let $\k$ be a field; the local ring $R=\poly{x,y}/(x^2,xy,y^3)$ has
  Hilbert series $1 + 2\tau + \tau^2$. The sequence
  \begin{equation*}
    R^2 \xra{(x\;y)} R \xra{x} R \xra{\binom{x}{y}} R^2
  \end{equation*}
  is exact and remains exact after dualization, but the product
  $(x,y)x$ is zero.
\end{exa}
\pagebreak
\begin{thm}
  \label{thm:ezdexist}
  Let $\Rm$ be a local ring with $\m^3=0$ and $e=\embdim R \ge 3$.
  Let $x$ be an element of $\mnotmm$; the following conditions are
  equivalent.
  \begin{eqc}
  \item The element $x$ is an exact zero divisor in $R$.
  \item The Hilbert series of $R$ is $1 + e\tau + (e-1)\tau^2$, and
    there exists an exact sequence of finitely generated free
    $R$-modules
    \begin{equation*}
      F:\quad F_3 \lra F_2 \lra F_1 \xra{\psi} F_0 \lra F_{-1}
    \end{equation*}
    such that $\Hom{F}{R}$ is exact, the homomorphisms are represented
    by matrices with entries in $\m$, and $\psi$ is represented by a
    matrix in which some row has $x$ as an entry and no other entry
    from $\mnotmm$.
  \end{eqc}
\end{thm}

\begin{prf*}
  If $x$ is an exact zero divisor in $R$, then the complex
  \eqref{epzdtac} supplies the desired exact sequence, and $R$ has
  Hilbert series $1 +e\tau +(e-1)\tau^2$; see \pgref{m3epzd}.

  To prove the converse, let $\Psi=(x_{ij})$ be a matrix of size $m
  \times n$ that represents $\psi$ and assume, without loss of
  generality, that the last row of $\Psi$ has exactly one entry
  $x=x_{mq}$ from $\mnotmm$.  By \partlemref{pm}{b} there is an
  equality $x\m = \m^2$ and, therefore, the length of $(x)$ is $e$
  by \partlemref{ann}{a}. As $R$ has length $2e$, additivity of length
  on short exact sequences yields $\lgt{\ann{x}} = e$.

  Let $(w_{ij})$ be an $n \times p$ matrix that represents the
  homomorphism $F_2 \to F_1$. The matrix equality $(x_{ij})(w_{ij})=0$
  yields $xw_{qj}=0$ for $j \in \set{1,\ldots,p}$; it follows that the
  ideal $\mathfrak{r} = (w_{q1},\ldots,w_{qp})$ is contained in
  $\ann{x}$. By \partlemref{pm}{b} some entry $w=w_{ql}$ is in
  $\mnotmm$, and there are inclusions
  \begin{equation}
    \label{eq:xr}
    (w) + \m^2 \subseteq \mathfrak{r} \subseteq \ann{x}.
  \end{equation}
  Now the inequalities
  \begin{equation*}
    e \le \lgt{((w)+\m^2)} \le \lgt{\ann{x}} = e
  \end{equation*}
  imply that equalities hold throughout \eqref{xr}; in particular,
  $(w) + \m^2 = \mathfrak{r}$ holds. This equality
  and \partlemref{pm}{b} yield $w\m = \m\mathfrak{r} = \m^2$; hence
  $w$ and $x$ form an exact pair of zero divisors
  by \partlemref{ann}{c}.
\end{prf*}

\begin{cor}
  \label{cor:diagonal}
  Let $\Rm$ be a local ring with $\m^3=0$. If $R$ is not Gorenstein,
  then the following conditions are equivalent:
  \begin{eqc}
  \item There is an exact zero divisor in $R$.
  \item For every $n \in \NN$ there is an indecomposable totally
    reflexive $R$-module that is presented by an upper triangular
    $n\times n$ matrix with entries in $\m$.
  \item There is a totally reflexive $R$-module without free summands
    that is presented by a matrix with entries in $\m$ and a
    row/column with only one entry in $\mnotmm$.
  \end{eqc}
\end{cor}

\begin{prf*}
  Set $e=\embdim R$; it is at least $2$ as $R$ is not Gorenstein. If
  $(i)$ holds, then $e$ is at least $3$, see \pgref{m3epzd}, so $(ii)$
  follows from \thmref{m3zd}. It is clear from \lemref{pm} that
  $(iii)$ follows from $(ii)$. To prove that $(iii)$ implies $(i)$,
  let $\Psi$ be a presentation matrix for a totally reflexive
  $R$-module without free summands and assume---possibly after
  replacing $\Psi$ with its transpose, which also presents a totally
  reflexive module---that some row of $\Psi$ has only one entry in
  $\mnotmm$.  By \pgref{m3tr} the Hilbert series of $R$ is $1 + e\tau
  + (e-1)\tau^2$, and $e$ is at least $3$, so it follows from
  \thmref{ezdexist} that there is an exact zero divisor in $R$.
\end{prf*}

The corollary manifests a strong relation between the existence of
exact zero divisors in $R$ and existence of totally reflexive
$R$-modules of any size.  A qualitatively different relation is
studied later; see~\rmkref{generic}.  These relations notwithstanding,
totally reflexive modules may exist in the absence of exact zero
divisors. In \secref{example} we exhibit a local ring $\Rm$, which has
no exact zero divisors, and a totally reflexive $R$-module that is
presented by a $2 \times 2$ matrix with all four entries from
$\mnotmm$. Thus the condition on the entries of the matrix
in \partcorref[]{diagonal}{$iii$} is~sharp.

\section{Families of non-isomorphic modules of the same size} 
\label{sec:bt2}

\noindent
In this section we continue the analysis of \conref{indec}.

\begin{dfn}
  Let $\Rmk$ be a local ring. Given a subset $\mathcal{K} \subseteq
  \k$, a subset of $R$ that contains exactly one lift of every element
  in $\mathcal{K}$ is called a \emph{lift} of $\mathcal{K}$ in $R$.
\end{dfn}

\begin{thm}
  \label{thm:noniso}
  Let $\Rmk$ be a local ring and let $\mathcal{L}$ be a lift of $\k$
  in $R$. Let $w$, $x$, $y$, $y'$, and $z$ be elements in $\mnotmm$
  and let $n$ be an integer. Assume that $w$ and $x$ form an exact
  pair of zero divisors in $R$ and that one of the following holds.
  \begin{prt}
  \item $n=2$ and the elements $w$, $x$, $y$, and $y'$ are linearly
    independent modulo $\m^2$;
  \item $n=2$, the elements $x$, $y$, and $y'$ are linearly
    independent modulo $\m^2$, and the element $w$ belongs to
    $(x)+\m^2$;
  \item $n \ge 3$, the elements $w$, $x$, $y$, and $y'$ are linearly
    independent modulo $\m^2$, and the following hold: $z\not\in (w) +
    \m^2$, $z\not\in (x)+ \m^2$, and $(y,y') \subseteq \ann{z}$;
    \text{ or}
  \item $n \ge 3$, the elements $x$, $y$, and $y'$ are linearly
    independent modulo $\m^2$, and the following hold: $w \in (x) +
    \m^2$, $z\not\in (x)+ \m^2$, and $(y,y') \subseteq \ann{z}$.
  \end{prt}
  Then the modules in the family $\{M_n(w,x,\lambda y +
  y',z)\}_{\lambda\in\mathcal L}$ are indecomposable, totally
  reflexive, and pairwise non-isomorphic.
\end{thm}

The proof of \thmref[]{noniso} takes up the balance of this section;
here is the cornerstone:

\begin{prp}
  \label{prp:noniso}
  Let $\Rm$ be a local ring and let $w$, $x$, $y$, $y'$, and $z$ be
  elements in $\mnotmm$. Assume that the following hold:
  \begin{equation*}
    z\not\in (w) + \m^2,\quad z\not\in (x)+ \m^2,\quad   y\not\in(w,x)+
    \m^2,  \qand y' \notin (w,y) + \m^2.
  \end{equation*}
  If $\theta$ and $\lambda$ are elements in $R$ with $\theta - \lambda
  \not\in \m$, and $n \ge 3$ is an integer, then the $R$-modules
  $M_n(w,x,\theta y + y', z)$ and $M_n(w,x,\lambda y + y', z)$ are
  non-isomorphic.
\end{prp}

The next example shows that the condition $n \ge 3$ in
\prpref[]{noniso} cannot be relaxed.

\begin{exa}
  Let $\Rm$ be a local ring with $\embdim R\ge 3$ and assume that $2$
  is a unit in $R$. Let $w$, $x$, and $y$ be linearly independent
  modulo $\m^2$, and set $y' = y - 2^{-1}x$. The equality
  \begin{equation*}
    \begin{pmatrix}
      1 & 1 \\ 0 & -1
    \end{pmatrix}
    \begin{pmatrix}
      w & y' \\ 0 & x
    \end{pmatrix}
    =
    \begin{pmatrix}
      w & y' -2y \\ 0 & x
    \end{pmatrix}
    \begin{pmatrix}
      1 & 0 \\ 0 & -1
    \end{pmatrix}
  \end{equation*}
  shows that the $R$-modules $M_2(w,x,0y+y')$ and $M_2(w,x,-2y+y')$
  are isomorphic.
\end{exa}

To get a statement similar to \prpref{noniso} for $2$-generated
modules, it suffices to assume that $y'$ is outside the span of $w$,
$x$, and $y$ modulo $\m^2$.
\pagebreak
\begin{prp}
  \label{prp:noniso2}
  Let $\Rm$ be a local ring and let $w$, $x$, $y$, and $y'$ be
  elements in $\mnotmm$. Assume that one of the following conditions
  holds:
  \begin{prt}
  \item The elements $w$, $x$, $y$, and $y'$ are linearly independent
    modulo $\m^2$.
  \item The elements $x$, $y$, and $y'$ are linearly independent
    modulo $\m^2$, and the element $w$ belongs to $(x)+\m^2$.
  \end{prt}
  If $\theta$ and $\lambda$ are elements in $R$ with $\theta - \lambda
  \not\in \m$, then the $R$-modules $M_2(w,x,\theta y + y')$ and
  $M_2(w,x,\lambda y + y')$ are non-isomorphic.
\end{prp}

\begin{prf*}
  Let $\theta$ and $\lambda$ be in $R$. Assume that (a) holds and that
  the $R$-modules $M_2(w,x,\theta y + y')$ and $M_2(w,x,\lambda y +
  y')$ are isomorphic. It follows that there exist matrices $\Alpha$
  and $\Beta$ in $\GL[R]{2}$ such that the equality
  \begin{equation*}
    \Alpha(w\Iota^o + x\Iota^e +(\theta y + y')\J^o) = (w\Iota^o +
    x\Iota^e +(\lambda y + y')\J^o)\Beta
  \end{equation*}
  holds; here the matrices $\Iota^o = \Iota_2^o$, $\Iota^e =
  \Iota_2^e$, and $\J^o = \J_2^o$ are as defined in
  \conref{indec}. The goal is to prove that $\theta - \lambda$ is in
  $\m$. After reduction modulo $\m$, the equality above yields, in
  particular,
  \begin{equation*}
    \bar{\Alpha}\J^o - \J^o\widebar{\Beta} =0=\bar{\theta}
    \bar{\Alpha}\J^o - \bar{\lambda} \J^o\widebar{\Beta},
  \end{equation*}
  and, therefore,
  $(\bar{\theta}-\bar{\lambda})\J^o\widebar{\Beta}=0$. As the matrix
  $\widebar{\Beta}$ is invertible and $\J^o$ is non-zero, this implies
  that $\theta - \lambda$ is in $\m$ as desired.

  If (b) holds, the desired conclusion is proved under Case 1 in the
  next proof.
\end{prf*}

\begin{prf*}[Proof of \pgref{prp:noniso}]
  Let $\theta$ and $\lambda$ be elements in $R$ and assume that
  $M_n(w,x,\theta y + y', z)$ and $M_n(w,x,\lambda y + y', z)$ are
  isomorphic as $R$-modules. It follows that there exist matrices
  $\Alpha$ and $\Beta$ in $\GL[R]{n}$ such that the equality
  \begin{equation}
    \label{eq:BT2}
    \Alpha(w\Iota^o + x\Iota^e +(\theta y + y')\J^o + z\J^e) =
    (w\Iota^o + x\Iota^e +(\lambda y + y')\J^o + z\J^e)\Beta
  \end{equation}
  holds; here the matrices $\Iota^o = \Iota_n^o$, $\Iota^e =
  \Iota_n^e$, $\J^o = \J_n^o$, and $\J^e = \J_n^e$ are as defined in
  \conref{indec}.  The goal is to prove that $\theta - \lambda$ is in
  $\m$.

  \subparagraph \emph{Case 1: $w$ is in $(x) + \m^2$.}  Under this
  assumption, one can write $w = rx + \delta$, where $\delta \in
  \m^2$, and rewrite \eqref{BT2} as
  \begin{equation}
    \label{eq:BT2a}
    \begin{split}
      x(r(\Alpha\Iota^o-\Iota^o\Beta) + \Alpha\Iota^e- \Iota^e\Beta)
      &+ y(\theta \Alpha\J^o - \lambda \J^o\Beta)\\
      &+ y'(\Alpha\J^o - \J^o\Beta) + z (\Alpha\J^e-\J^e\Beta) =
      \Delta,
    \end{split}
  \end{equation}
  where $\Delta$ is a matrix with entries in $\m^2$. The assumptions
  on $w$, $x$, $y$, and $y'$ imply $y \notin (x) + \m^2$ and $y'
  \notin (x,y) + \m^2$, so the elements $x$, $y$, and $y'$ are
  linearly independent modulo $\m^2$. There exist elements $v_i$ such
  that $v_1,\ldots,v_{e-3},x,y,y'$ form a minimal set of generators
  for $\m$. Write
  \begin{equation*}
    z = sx + ty + uy' + \sum_{i=1}^{e-3}d_iv_i,
  \end{equation*}
  substitute this expression into \eqref{BT2a}, and reduce modulo $\m$
  to get
  \begin{align}
    \label{eq:rs}
    0 &= \bar{r}(\bar{\Alpha}\Iota^o - \Iota^o\widebar{\Beta}) +
    \bar{\Alpha}\Iota^e - \Iota^e\widebar{\Beta} +
    \bar{s}(\bar{\Alpha}\J^e-\J^e\widebar{\Beta}),\\
    \label{eq:-t}
    0 &= \bar{\theta} \bar{\Alpha}\J^o - \bar{\lambda}
    \J^o\widebar{\Beta}
    + \bar{t}(\bar{\Alpha}\J^e-\J^e\widebar{\Beta}),\\
    \label{eq:-u}
    0 &= \bar{\Alpha}\J^o - \J^o\widebar{\Beta}
    + \bar{u}(\bar{\Alpha}\J^e-\J^e\widebar{\Beta}),\text{ and}\\
    \label{eq:-d}
    0 &= \bar{d}_i(\bar{\Alpha}\J^e-\J^e\widebar{\Beta}), \text{ for
      $i \in \set{1,\ldots,e-3}$.}
  \end{align}
  The arguments that follow use the relations that
  \eqref{rs}--\eqref{-d} induce between the entries of the matrices
  $\bar{\Alpha} = (a_{ij})$ and $\widebar{\Beta} = (b_{ij})$. With the
  convention $a_{hl}=0=b_{hl}$ for $h,l \in \set{0,n+1}$, it is
  elementary to verify that~the following systems of equalities hold
  for $i$ and $j$ in $\set{1,\ldots,n}$ and elements $f$ and $g$
  in~$R/\m$:
  \begin{align}
    \label{eq:BT2Io}
    (\bar{\Alpha}\Iota^o - \Iota^o\widebar{\Beta})_{ij} &=
    \begin{cases}
      a_{ij}-b_{ij} & \text{$i$ odd, $j$ odd}\\
      -b_{ij} & \text{$i$ odd, $j$ even}\\
      a_{ij} & \text{$i$ even, $j$ odd}\\
      0 & \text{$i$ even, $j$ even}
    \end{cases}\\
    \label{eq:BT2Ie}
    (\bar{\Alpha}\Iota^e - \Iota^e\widebar{\Beta})_{ij} &=
    \begin{cases}
      0 & \text{$i$ odd, $j$ odd}\\
      a_{ij} & \text{$i$ odd, $j$ even}\\
      -b_{ij} & \text{$i$ even, $j$ odd}\\
      a_{ij}-b_{ij} & \text{$i$ even, $j$ even}
    \end{cases}\\
    \label{eq:BT2Jo}
    (f\bar{\Alpha}\J^o - g\J^o\widebar{\Beta})_{ij} &=
    \begin{cases}
      gb_{(i+1)j} & \text{$i$ odd, $j$ odd}\\
      fa_{i(j-1)} - gb_{(i+1)j} & \text{$i$ odd, $j$ even}\\
      0 & \text{$i$ even, $j$ odd}\\
      fa_{i(j-1)} & \text{$i$ even, $j$ even}
    \end{cases}\\
    \label{eq:BT2Je}
    (\bar{\Alpha}\J^e - \J^e\widebar{\Beta})_{ij} &=
    \begin{cases}
      a_{i(j-1)} & \text{$i$ odd, $j$ odd}\\
      0  & \text{$i$ odd, $j$ even}\\
      a_{i(j-1)} - b_{(i+1)j}  & \text{$i$ even, $j$ odd}\\
      -b_{(i+1)j} & \text{$i$ even, $j$ even}
    \end{cases}
  \end{align}
  Each of the equalities \eqref{rs}--\eqref{-d} induces four
  subsystems, which are referred to by subscripts `oo', `oe', `eo',
  and `ee', where `oo' stands for '$i$ odd and $j$ odd' etc. For
  example, \subeqref{-t}{oe} refers to the equalities $0 =
  \bar{\theta}a_{i(j-1)} - \bar{\lambda}b_{(i+1)j}$, for $i$ odd and
  $j$ even.

  Let $n \ge 2$; the goal is to prove the equality
  $\bar{\theta}=\bar{\lambda}$, as that implies $\theta - \lambda \in
  \m$. First assume $\bar{d}_i\ne 0$ for some $i \in
  \set{1,\ldots,e-3}$, then \eqref{-d} yields
  $\bar{\Alpha}\J^e-\J^e\widebar{\Beta}=0$. From \eqref{-t} and
  \eqref{-u} one then gets
  \begin{equation*}
    \bar{\Alpha}\J^o - \J^o\widebar{\Beta} =0=\bar{\theta}
    \bar{\Alpha}\J^o - \bar{\lambda} \J^o\widebar{\Beta}
  \end{equation*}
  and thus $(\bar{\theta}-\bar{\lambda})\J^o\widebar{\Beta}=0$. As
  $\widebar{\Beta}$ is invertible and $\J^o$ is non-zero, this yields
  the desired equality $\bar{\theta}=\bar{\lambda}$. Henceforth we
  assume $\bar{d}_i = 0$ for all $i \in \set{1\ldots,e-3}$.

  By assumption, $z$ is not in $(x)+\m^2$, so $\bar{u}$ or $\bar{t}$
  is non-zero. In case $\bar{u}$ is non-zero, \subeqref{-u}{oo} and
  \subeqref{-u}{eo} yield $b_{h1}=0$ for $1 < h \le n$. If $\bar{t}$
  is non-zero, then \subeqref{-t}{oo} and \subeqref{-t}{eo} yield the
  same conclusion. The matrix $\widebar{\Beta}$ is invertible, so each
  of its columns contains a non-zero element. It follows that $b_{11}$
  is non-zero, and by \subeqref{rs}{oo} one has $a_{11}=b_{11}$. From
  \subeqref{-t}{oe} and \subeqref{-u}{oe} one gets $a_{11} = b_{22}$
  and $(\bar{\theta} - \bar{\lambda})a_{11} = 0$, whence the equality
  $\bar{\theta}=\bar{\lambda}$ holds.

  For $n=2$ the arguments above establish the assertion in
  \prpref{noniso2} under assumption (b) \emph{ibid.}

\enlargethispage*{\baselineskip}
  \subparagraph\emph{Case 2: $w$ is not in $(x)+\m^2$.} It follows
  from the assumption $y \notin (w,x) + \m^2$ that $w$, $x$, and $y$
  are linearly independent modulo $\m^2$, so there exist elements
  $v_i$ such that $v_1,\ldots,v_{e-3},w,x,y$ form a minimal set of
  generators for $\m$. Write
  \begin{equation*}
    y' = pw + qx + ry + \sum_{i=1}^{e-3}c_iv_i \qqand z = sw + tx + uy
    + \sum_{i=1}^{e-3}d_iv_i.
  \end{equation*}
  Substitute these expressions into \eqref{BT2} and reduce modulo $\m$
  to get
  \begin{align}
    \label{eq:ps}
    0 &= \bar{\Alpha}\Iota^o - \Iota^o\widebar{\Beta} +
    \bar{p}(\bar{\Alpha}\J^o - \J^o\widebar{\Beta}) +
    \bar{s}(\bar{\Alpha}\J^e-\J^e\widebar{\Beta}),\\
    \label{eq:qt}
    0 &=\bar{\Alpha}\Iota^e - \Iota^e\widebar{\Beta} +
    \bar{q}(\bar{\Alpha}\J^o - \J^o\widebar{\Beta}) +
    \bar{t}(\bar{\Alpha}\J^e-\J^e\widebar{\Beta}),\\
    \label{eq:ru}
    0 &= \bar{\theta} \bar{\Alpha}\J^o - \bar{\lambda}
    \J^o\widebar{\Beta} + \bar{r}(\bar{\Alpha}\J^o -
    \J^o\widebar{\Beta}) +
    \bar{u}(\bar{\Alpha}\J^e-\J^e\widebar{\Beta}),\text{ and}\\
    \label{eq:cd}
    0 &= \bar{c_i}(\bar{\Alpha}\J^o - \J^o\widebar{\Beta}) +
    \bar{d}_i(\bar{\Alpha}\J^e-\J^e\widebar{\Beta}), \text{ for $i \in
      \set{1,\ldots,e-3}$}.
  \end{align}
  As in Case 1, set $\bar{\Alpha} = (a_{ij})$ and $\widebar{\Beta} =
  (b_{ij})$ so that the equalities \eqref{BT2Io}--\eqref{BT2Je} hold.
  The subscripts `oo', `oe', `eo', and `ee' are used, as in Case 1, to
  denote the subsystems induced by \eqref{ps}--\eqref{cd}.

  Let $n\ge 3$; the goal is, again, to prove the equality
  $\bar{\theta}=\bar{\lambda}$. In the following, $h$ and $l$ are
  integers in $\set{1,\ldots,n}$. First, notice that if $\bar{c}_m \ne
  0$ for some $m \in \set{1,\ldots,e-3}$, then \subeqref{cd}{oe}
  implies $a_{11} = b_{22}$ and, in turn, \subeqref{ru}{oe} yields
  $(\bar{\theta} - \bar{\lambda})a_{11}=0$. To conclude
  $\bar{\theta}=\bar{\lambda}$, it must be verified that $a_{11}$ is
  non-zero.  To this end, assume first $\bar{d}_m=0$; from
  \subeqref{cd}{ee} and \subeqref{cd}{oo} one immediately gets
  \begin{equation}
    \label{eq:abk1}
    a_{h1} = 0 = b_{h1} \quad\text{for $h$ even.}
  \end{equation}
  As $z$ is in $\m\backslash\m^2$, one of the coefficients
  $\bar{d}_1,\ldots,\bar{d}_{e-3},\bar{s},\bar{t},\bar{u}$ is
  non-zero. If $\bar{u}$ or one of $\bar{d}_1,\ldots,\bar{d}_{e-3}$ is
  non-zero, then \subeqref{ru}{eo} or \subeqref{cd}{eo} yields
  $b_{h1}=0$ for $h>1$ odd. If $\bar{s}$ or $\bar{t}$ is non-zero,
  then the same conclusion follows from \eqref{abk1} combined with
  \subeqref{ps}{eo} or with \subeqref{qt}{eo}. Now that $b_{h1}=0$
  holds for all $h \ge 2$, it follows that $b_{11}$ is
  non-zero. Finally, \subeqref{ps}{oo} yields $a_{11}=b_{11}$, so the
  entry $a_{11}$ is not zero, as desired. Now assume $\bar{d}_m\ne 0$,
  then \subeqref{cd}{oo} and \subeqref{cd}{eo} immediately give
  $b_{h1} = 0$ for $h>1$.  As above, we conclude that $b_{11}$ is
  non-zero, whence $a_{11} \ne 0$ by \subeqref{ps}{oo}.  This
  concludes the argument under the assumption that one of the
  coefficients $\bar{c}_i$ is non-zero. Henceforth we assume
  $\bar{c}_i = 0$ for all $i \in \set{1,\ldots,e-3}$; it follows that
  $\bar{q}$ is non-zero, as $y' \notin (w,y) + \m^2$ by assumption.

  If $\bar{d}_i \ne 0$ for some $i \in \set{1,\ldots,e-3}$, then
  \subeqref{cd}{oo} yields $a_{12}=0$, as $n$ is at least $3$. From
  \subeqref{qt}{oe} one gets $a_{11}=b_{22}$, and then
  \subeqref{ru}{oe} implies $(\bar{\theta} -
  \bar{\lambda})a_{11}=0$. To see that $a_{11}$ is non-zero, notice
  that \subeqref{cd}{eo} yields $b_{h1}=0$ for $h>1$ odd, while
  \subeqref{qt}{oo} yields $b_{h1}=0$ for $h$ even. It follows that
  $b_{11}$ is non-zero, and then $a_{11}\ne 0$, by
  \subeqref{ps}{oo}. Thus the desired equality $\bar{\theta} =
  \bar{\lambda}$ holds.  This concludes the argument under the
  assumption that one of the coefficients $\bar{d}_i$ is
  non-zero. Henceforth we assume $\bar{d}_i = 0$ for all $i \in
  \set{1,\ldots,e-3}$.

  To finish the argument, we deal separately with the cases $\bar{u}
  \ne 0$ and $\bar{u} = 0$. Assume first $\bar{u} \ne 0$. By
  assumption one has $n \ge 3$, so \subeqref{ru}{eo} yields
  $a_{22}=b_{33}$, and then it follows from \subeqref{qt}{eo} that
  $b_{23}$ is zero. From \subeqref{ru}{oo} one gets $a_{12}=0$, and
  then the equality $a_{11} = b_{22}$ follows from
  \subeqref{qt}{oe}. In turn, \subeqref{ru}{oe} implies
  $(\bar{\theta}- \bar{\lambda})a_{11} = 0$. To see that $a_{11}$ is
  non-zero, notice that there are equalities $b_{h1}=0$ for $h>1$ odd
  by \subeqref{ru}{eo} and $b_{h1}=0$ for $h$ even by
  \subeqref{qt}{oo}. It follows that $b_{11}$ is non-zero, and the
  equality $a_{11}=b_{11}$ holds by \subeqref{ps}{oo}. As above we
  conclude $\bar{\theta} = \bar{\lambda}$.

  Finally, assume $\bar{u} = 0$; the assumptions $z\notin (w) + \m^2$
  and $z\notin (x) + \m^2$ imply that $\bar{s}$ and $\bar{t}$ are both
  non-zero. From \subeqref{ru}{ee} and \subeqref{ru}{oo} one gets:
  \begin{equation}
    \label{eq:abr}
    (\bar{\theta} + \bar{r})a_{h1} = 0 = (\bar{\lambda} +
    \bar{r})b_{h1} \quad\text{for $h$ even.}
  \end{equation}
  First assume that $\bar{\lambda} + \bar{r}$ is zero. If
  $\bar{\theta} + \bar{r}$ is zero, then the desired equality
  $\bar{\theta} = \bar{\lambda}$ holds. If $\bar{\theta} + \bar{r}$ is
  non-zero, then \eqref{abr} gives $a_{h1}=0$ for $h$ even. Moreover,
  \subeqref{ru}{oe} implies $(\bar{\theta} + \bar{r})a_{h1} =
  (\bar{\lambda} + \bar{r})b_{(h+1)2} = 0$, so $a_{h1}=0$ for $h$ odd
  as well, which is absurd as $\bar{\Alpha}$ is invertible. Now assume
  that $\bar{\lambda} + \bar{r}$ is non-zero. From \eqref{abr} one
  gets $b_{h1}=0$ for $h$ even and, in turn, \subeqref{qt}{eo} yields
  $b_{h1}=0$ for $h>1$ odd. It follows that $b_{11}$ is non-zero, and
  then one has $a_{11}\ne 0$ by \subeqref{ps}{oo}. By assumption, $n$
  is at least $3$, so \subeqref{ru}{oo} gives $(\bar{\lambda} +
  \bar{r})b_{23}=0$, which implies $b_{23}=0$. Now \subeqref{qt}{oo}
  yields $a_{12}=0$, and then \subeqref{qt}{oe} implies
  $a_{11}=b_{22}$. From \subeqref{ru}{oe} one gets $(\bar{\theta} -
  \bar{\lambda})a_{11}$, whence the equality $\bar{\theta} =
  \bar{\lambda}$ holds.
\end{prf*}

\begin{prf*}[Proof of \pgref{thm:noniso}]
  Under the assumptions in part (a) or (b), it is immediate from
  \thmref{indectr} that the module $M_2(w,x,\lambda y + y')$ is
  indecomposable and totally reflexive for every $\lambda \in R$. It
  follows from \prpref{noniso2} that the modules in the family
  $\{M_2(w,x,\lambda y + y')\}_{\lambda\in\mathcal L}$ are pairwise
  non-isomorphic.

  Fix $n\ge 3$. \prpref{noniso} shows that, under the assumptions in
  part (c) or (d), the modules in the family $\{M_n(w,x,\lambda
  y+y',z)\}_{\lambda\in\mathcal L}$ are pairwise
  non-isomorphic. Moreover, for every $\lambda \in R$, one has
  $(\lambda y +y')z=0$, so it follows from \thmref{indectr} the module
  $M_n(w,z,\lambda y + y',z)$ is totally reflexive and indecomposable.
\end{prf*}

\section{Brauer--Thrall II over short local rings with exact zero
  divisors}
\label{sec:BT2}

\noindent
In this section $\Rmk$ is a local ring with $\m^3=0$. Together,
\eqref{length} and the theorems \thmref[]{indectr}, \thmref[]{bt2},
\thmref[]{BT21}, and \thmref[]{BT22} establish \thmref{secondbt}.

\begin{rmk}
  \label{rmk:2gor}
  Assume that $R$ has embedding dimension $2$. If there is an exact
  zero divisor in $R$, then the Hilbert series of $R$ is $1 + 2\tau +
  \tau^2$, and the equality $\ann{\m} = \m^2$ holds; see
  \pgref{m3epzd}. Therefore, $R$ is Gorenstein, and the equality $x\m=
  \m^2$ holds for all $x \in \mnotmm$; in particular, every element in
  $\m\backslash \m^2$ is an exact zero divisor
  by \partlemref{ann}{c}. On the other hand, if $R$ is Gorenstein,
  then one has $\operatorname{H}_R(\tau) = 1 + 2\tau + \tau^2$ and
  $\ann{\m} = \m^2$. It is now elementary to verify that the
  equalities $\lgt{(x)} = 2 = \lgt{\ann{x}}$ hold for every $x\in
  \mnotmm$, so every such element is an exact zero divisor in $R$
  by \partlemref{ann}{c}. It follows from work of Serre
  \prpcite[5]{JPS60} that $R$ is Gorenstein if and only if it is
  complete intersection; thus the following conditions are equivalent:
  \begin{eqc}
  \item There is an exact zero divisor in $R$.
  \item Every element in $\mnotmm$ is an exact zero divisor.
  \item $R$ is complete intersection.
  \end{eqc}
\end{rmk}

In contrast, if $R$ has embedding dimension at least $3$, and $\k$ is
algebraically closed, then there exist elements in $\mnotmm$ that are
not exact zero divisors. This fact follows from \lemref{nezd}, and it
is essential for our proof of \thmref{bt2}.

\begin{rmk}
  \label{rmk:rank}
  Assume that $R$ has Hilbert series $1 +e\tau +f\tau^2$, and let $x$
  be an element in $\m$. The equality $x\m = \m^2$ holds if and only
  if the $\k$-linear map from $\m/\m^2$ to $\m^2$ given by
  multiplication by $x$ is surjective. Let $\Xi_x$ be a matrix that
  represents this map; it is an $f \times e$ matrix with entries in
  $\k$, so the equality $x\m = \m^2$ holds if and only if $\Xi_x$ has
  rank $f$.

  Assume that the (in)equalities $f = e-1 \ge 1$ hold. If $w$ and $x$
  are elements in $R$ with $wx=0$, then they form an exact pair of
  zero divisors if and only if both matrices $\Xi_x$ and $\Xi_w$ have
  a non-zero maximal minor; cf.~\partlemref{ann}{c}.
\end{rmk}

\begin{lem}
  \label{lem:nezd}
  Assume that $\Rmk$ has Hilbert series $1 +e\tau +f\tau^2$ and that
  $\k$ is algebraically closed; set $n=e-f+1$. If one has $2\le f \le
  e-1$ and $v_0,\ldots,v_n\in \m$ are linearly independent modulo
  $\m^2$, then there exist $r_0,\ldots,r_n\in R$, at least one of
  which is invertible, such that the ideal $(\sum_{h=0}^nr_hv_h)\m$ is
  properly contained in $\m^2$.
\end{lem}

\begin{prf*}
  Let $\set{x_1,\ldots,x_e}$ be a minimal set of generators for $\m$
  and let $\set{u_1,\ldots,u_f}$ be a basis for the $\k$-vector space
  $\m^2$.  For $h\in\set{0,\ldots, n}$ and $j\in\set{1,\ldots, e}$
  write
  \begin{equation*}
    v_hx_j = \sum_{i=1}^{f} \xi_{hij} u_i \,,
  \end{equation*}
  where the elements $\xi_{hij}$ are in $\k$.  In the following,
  $\bar{r}$ denotes the image in $\k$ of the element $r\in R$. For
  $r_0,\ldots,r_n$ in $R$ and $v = \sum_{h=0}^n r_hv_h$ the equality
  $v\m = \m^2$ holds if and only if the $f \times e$ matrix
  \begin{equation*}
    \Xi_v = \left( \sum_{h=0}^n \bar{r}_h\xi_{hij} \right)_{ij}
  \end{equation*}
  has rank $f$; see \rmkref{rank}. Let $\Chi$ be the matrix obtained
  from $\Xi_v$ by replacing $\bar{r}_0,\ldots,\bar{r}_n$ with
  indeterminates $\chi_0,\ldots,\chi_n$. The non-zero entries in
  $\Chi$ are then homogeneous linear forms in $\chi_0,\ldots,\chi_n$,
  and in the polynomial algebra $\poly[\k]{\chi_0,\ldots,\chi_n}$ the
  ideal $I_{f}(X)$, generated by the maximal minors of $\Chi$, has
  height at most $n$; see \thmcite[(2.1)]{bruvet}. As $I_{f}(\Chi)$ is
  generated by homogeneous polynomials, it follows from Hilbert's
  Nullstellensatz that there exists a point $(\rho_0, \ldots, \rho_n)$
  in $\k^{n+1}\backslash \{0\}$ such that all maximal minors of the
  matrix $\smash{\left( \sum_{h=0}^n \rho_h\xi_{hij} \right)_{ij}}$
  vanish. Let $r_0, \ldots, r_n$ be lifts of $\rho_0, \ldots, \rho_n$
  in $R$, then at least one of them is not in $\m$, and the ideal
  $(\smash\sum_{h=0}^nr_hv_h)\m$ is properly contained in $\m^2$.
\end{prf*}

\begin{thm}
  \label{thm:bt2}
  Let $\Rmk$ be a local ring with $\m^3=0$, $\embdim R \ge 3$, and
  $\k$ algebraically closed. If $w$ and $x$ form an exact pair of zero
  divisors in $R$, then there exist elements $y$, $y'$, and $z$ in
  $\mnotmm$, such that for every lift $\mathcal{L}$ of
  $\k\backslash\{0\}$ in $R$ and for every integer $n \ge 3$ the
  modules in the family $\{M_n(w,x,\lambda
  y+y',z)\}_{\lambda\in\mathcal L}$ are indecomposable, totally
  reflexive, and pairwise non-isomorphic.
\end{thm}

\begin{prf*}
  Assume that $w$ and $x$ form an exact pair of zero divisors in
  $R$. By \pgref{m3epzd} the Hilbert series of $R$ is $1 + e\tau +
  (e-1)\tau^2$, so it follows from \lemref{nezd} that there exists an
  element $z \in \mnotmm$ such that $z\m$ is properly contained in
  $\m^2$. In particular, one has $\lgt{(z)} < e$ and, therefore,
  $\lgt{\ann{z}} > e$ by additivity of length on short exact
  sequences. It follows that there exist two elements, call them $y$
  and $y'$, in $\ann{z}$ that are linearly independent modulo
  $\m^2$. The inequality $\lgt{(z)}<e$ implies that $z$ is not in $(w)
  = (w) + \m^2$ and not in $(x) = (x) + \m^2$.  If $y$ and $y'$ were
  both in $(w,x) = (w,x) + \m^2$, then $x$ would be in $(y,y')$, which
  is impossible as $z \notin (w)$. Without loss of generality, assume
  $y\notin(w,x) + \m^2$. If $y'$ were in $(w,y) = (w,y) + \m^2$, then
  $w$ would be in $(y,y')$, which is also impossible. Now let
  $\mathcal{L}$ be a lift of $\k\backslash\set{0}$ in $R$ and let $n
  \ge 3$ be an integer. It follows from \prpref{noniso} that the
  modules in the family $\{M_n(w,x,\lambda
  y+y',z)\}_{\lambda\in\mathcal L}$ are pairwise non-isomorphic.

  If $w$ is in $(x) = (x) + \m^2$, then it follows from
  \thmref{indectr} that the modules in the family $\{M_n(w,x,\lambda
  y+y',z)\}_{\lambda\in\mathcal L}$ are indecomposable and totally
  reflexive. Indeed, for every $\lambda \in R$ the element $\lambda y
  + y'$ annihilates $z$, whence it is not in $(x) = (x) + \m^2$.

  If $w$ is not in $(x) = (x) + \m^2$, then it follows from the
  assumption $y\notin(w,x) + \m^2$ that $w$, $x$, and $y$ are linearly
  independent modulo $\m^2$. There exist elements $v_i$ such that
  $v_1,\ldots,v_{e-3},w,x,y$ form a minimal set of generators for
  $\m$. Write
  \begin{equation*}
    y' = pw + qx + ry + \sum_{i=1}^{e-3}c_iv_i.
  \end{equation*}
  As $yz=0$ holds and the elements $y$ and $y'$ are linearly
  independent modulo $\m^2$, we may assume $r=0$. For $\lambda \in
  \mathcal{L}$, the elements $w$, $x$, and $\lambda y + y'$ are
  linearly independent. Indeed, if there is a relation
  \begin{equation*}
    sw + tx + u(\lambda y + y') = (s+up)w + (t + uq)x + u\lambda y +
    u\sum_{i=1}^{e-3}c_iv_i \in \m^2,
  \end{equation*}
  then $u$ is in $\m$ as $\lambda\notin \m$, and then $s$ and $t$ are
  in $\m$ as $w$ and $x$ are linearly independent modulo $\m^2$. Now
  it follows from \thmref{indectr} that the modules in the family
  $\{M_n(w,x,\lambda y+y',z)\}_{\lambda\in\mathcal L}$ are
  indecomposable and totally reflexive.
\end{prf*}

\begin{rmk}
  Assume that $R$ has embedding dimension $3$ and that $\m$ is
  minimally generated by elements $v$, $w$, and $x$, where $w$ and $x$
  form an exact pair of zero divisors. Let $y$ and $y'$ be elements in
  $\m$ and let $\mathcal{L}$ be a subset of $R$. It is elementary to
  verify that each module in the family $\{M_2(w,x,\lambda
  y+y')\}_{\lambda\in\mathcal L}$ is isomorphic to either $R/(w)
  \oplus R/(x)$ or to the indecomposable $R$-module $M_2(w,x,v)$. Thus
  the requirement $n \ge 3$ in \thmref{bt2} cannot be relaxed.
\end{rmk}

We now proceed to deal with $1$- and $2$-generated modules.

\begin{thm}
  \label{thm:BT21}
  Let $\Rmk$ be a local ring with $\m^3=0$, $\embdim R\ge 2$, and $\k$
  infinite.  If there is an exact zero divisor in $R$, then the set
  \begin{equation*}
    \mathcal{N}_1 =
    \setof{R/(x)}{\text{$x$ is an exact zero divisor in $R$}}
  \end{equation*}
  of totally reflexive $R$-modules contains a subset $\mathcal{M}_1$
  of cardinality $\card{\k}$, such that the modules in $\mathcal{M}_1$
  are pairwise non-isomorphic.
\end{thm}

\begin{prf*}
  Assume that there is an exact zero divisor in $R$; then $R$ has
  Hilbert series $1+e\tau+(e-1)\tau^2$; see~\pgref{m3epzd}. Let
  $\set{x_1,\ldots,x_e}$ be a minimal set of generators for $\m$ and
  let $\set{u_1,\ldots,u_{e-1}}$ be a basis for $\m^2$.  For $h$ and
  $j$ in $\set{1,\ldots, e}$ write
  \begin{equation*}
    x_hx_j = \sum_{i=1}^{e-1} \xi_{hij} u_i,
  \end{equation*}
  where the elements $\xi_{hij}$ are in $\k$.  For $r_1,\ldots,r_e$ in
  $R$ let $\bar{r}_h$ denote the image of $r_h$ in $\k$, and set $x =
  r_1x_1+ \dots + r_ex_e$. The equality $x\m = \m^2$ holds if and only
  if the $(e-1) \times e$ matrix
  \begin{equation*}
    \Xi_x =  \left( \sum_{h=1}^e \bar{r}_h\xi_{hij} \right)_{ij}
  \end{equation*}
  has rank $e-1$; see \rmkref{rank}. For $j\in\set{1,\ldots,e}$ let
  $\mu_j(x)$ be the maximal minor of $\Xi_x$ obtained by omitting the
  $j$th column. Notice that each minor $\mu_j(x)$ is a homogeneous
  polynomial expression in the elements $\bar{r}_1,\ldots,
  \bar{r}_e$. The column vector
  $(\mu_1(x)\;-\!\mu_2(x)\;\cdots\;(-1)^{e-1}\mu_e(x))^\mathrm{T}$ is
  in the null-space of the matrix $\Xi_x$, so the element
  \begin{equation*}
    w=\mu_1(x)x_1-\mu_2(x)x_2+\cdots+(-1)^{e-1}\mu_e(x)x_e
  \end{equation*}
  annihilates $x$. Let $\nu_1(w),\ldots,\nu_e(w)$ be the maximal
  minors of the matrix $\Xi_w$; they are homogeneous polynomial
  expressions in $\mu_1(x),\ldots,\mu_e(x)$ and, therefore, in
  $\bar{r}_1,\ldots, \bar{r}_e$. By \rmkref{rank} the elements $x$ and
  $w$ form an exact pair of zero divisors if and only if both matrices
  $\Xi_x$ and $\Xi_w$ have rank $e-1$.

  For $j \in \set{1,\ldots,e}$ define $\mu_j$ and $\nu_j$ to be the
  homogeneous polynomials in $\poly[\k]{\chi_1,\ldots,\chi_e}$
  obtained from $\mu_j(x)$ and $\nu_j(w)$ by replacing the elements
  $\bar{r}_h$ by the indeterminates $\chi_h$, for $h \in
  \set{1,\ldots,e}$. In the projective space $\mathbb{P}_\k^{e-1}$,
  the complement $\mathcal{E}$ of the intersection of vanishing sets
  $\mathcal{Z}(\nu_1) \cap \cdots \cap \mathcal{Z}(\nu_e)$ is an open
  set. No point in $\mathcal{E}$ is in the intersection
  $\mathcal{Z}(\mu_1) \cap \cdots \cap \mathcal{Z}(\mu_e)$, as each
  polynomial $\nu_j$ is a polynomial in
  $\mu_1,\ldots,\mu_e$. Therefore, each point $(\rho_1: \cdots:
  \rho_e)$ in $\mathcal{E}$ corresponds to an exact zero divisor as
  follows. Let $r_1,\ldots,r_e$ be lifts of $\rho_1,\ldots,\rho_e$ in
  $R$; then the element $x=r_1x_1+\cdots+r_ex_e$ is an exact zero
  divisor. It is clear that two distinct points in $\mathcal{E}$
  correspond to non-isomorphic modules in $\mathcal{N}_1$. Take as
  $\mathcal{M}_1$ any subset of $\mathcal{N}_1$ such that the elements
  of $\mathcal{M}_1$ are in one-to-one correspondence with the points
  in $\mathcal{E}$.  By assumption, the set $\mathcal{E}$ is
  non-empty, and, if $\k$ is infinite, then $\mathcal{E}$ has the same
  cardinality as $\k$.
\end{prf*}

\begin{rmk}
  The modules in $\mathcal{M}_1$ are in one-to-one correspondence with
  the points of a non-empty Zariski open set in
  $\mathbb{P}_\k^{e-1}$. See also \rmkref{zar}.
\end{rmk}

\begin{thm}
  \label{thm:BT22}
  Let $\Rmk$ be a local ring with $\m^3=0$, $\embdim R \ge 3$, and
  $\k$ infinite. If there is an exact zero divisor in $R$, then the
  set
  \begin{equation*}
    \mathcal{N}_2 =
    \left\{\: M_2(w,x,y)
      \:\left|\:
        \begin{gathered}
          \text{ $w$, $x$, and $y$ are elements in
            $R$, such that} \\[-.6ex]
          \text{$w$ and $x$ form an exact pair of zero divisors }
        \end{gathered}
      \right.
    \right\}
  \end{equation*}
  of totally reflexive $R$-modules contains a subset $\mathcal{M}_2$
  of cardinality $\card{\k}$, such that the modules in $\mathcal{M}_2$
  are indecomposable and pairwise non-isomorphic.
\end{thm}

\begin{prf*}
  Assume that there is an exact zero divisor in $R$ and let
  $\mathcal{M}_1$ be the set of cyclic totally reflexive $R$-modules
  afforded by \thmref{BT21}; the cardinality of $\mathcal{M}_1$ is
  $\card{\k}$. From $\mathcal{M}_1$ one can construct another set of
  the same cardinality, whose elements are exact pairs of zero
  divisors, such that for any two of them, say $w,x$ and $w',x'$, one
  has $(x) \not\is (x')$ and $(w) \not\is (x')$. Given two such pairs,
  choose elements $y$ and $y'$ in $\mnotmm$ such that $y \notin (w,x)$
  and $y' \notin (w',x')$. By \thmref{m3zd}, the modules $M_2(w,x,y)$
  and $M_2(w',x',y')$ are indecomposable and totally reflexive.

  Suppose that the $R$-modules $M_2(w,x,y)$ and $M_2(w',x',y')$ are
  isomorphic; then there exist matrices $\Alpha = (a_{ij})$ and $\Beta
  = (b_{ij})$ in $\GL[R]{2}$ such that the equality
  \begin{equation*}
    \Alpha
    \begin{pmatrix}
      w & y\\ 0 & x
    \end{pmatrix}
    =
    \begin{pmatrix}
      w' & y'\\ 0 & x'
    \end{pmatrix}
    \Beta
  \end{equation*}
  holds. In particular, there are equalities
  \begin{equation*}
    a_{21}w = b_{21}x' \qqand a_{21}y + a_{22}x = b_{22}x'.
  \end{equation*}
  The first one shows that the entries $a_{21}$ and $b_{21}$ are
  elements in $\m$, and then it follows from the second one that
  $a_{22}$ and $b_{22}$ are in $\m$. Thus $\Alpha$ and $\Beta$ each
  have a row with entries in $\m$, which contradicts the assumption
  that they are invertible.
\end{prf*}

\section{Existence of exact zero divisors}
\label{sec:quad}

\noindent
In previous sections we constructed families of totally reflexive
modules starting from an exact pair of zero divisors. Now we address
the question of existence of exact zero divisors; in particular, we
prove \thmref{generic}; see \rmkref{generic}.

A local ring $\Rm$ with $\m^3=0$ and embedding dimension $1$ is, by
Cohen's Structure Theorem, isomorphic to $D/(d^2)$ or $D/(d^3)$, where
$(D,(d))$ is a discrete valuation domain. In either case, $d$ is an
exact zero divisor. In the following we focus on rings of embedding
dimension at least $2$.

\begin{rmk}
  Let $\Rm$ be a local ring with $\m^3=0$. It is elementary to verify
  that elements in $\m$ annihilate each other if and only if their
  images in the associated graded ring $\operatorname{gr}_\m(R)$
  annihilate each other. Thus an element $x\in\m$ is an exact zero
  divisor in $R$ if and only if $\bar{x} \in \m/\m^2$ is an exact zero
  divisor in $\operatorname{gr}_\m(R)$.
\end{rmk}

Let $\Rmk$ be a standard graded $\k$-algebra with $\m^3=0$ and
embedding dimension $e \ge 2$.  Assume that there is an exact zero
divisor in $R$; by \pgref{m3epzd} one has $\operatorname{H}_R(\tau) =
1 + e\tau + (e-1)\tau^2$. If $e$ is $2$, then it follows from
\cite[Hilfssatz~7]{GSc64} that $R$ is complete intersection with
Poincar\'e series
\begin{equation*}
  \sum_{i=0}^\infty \bet{i}{\k}\tau^i = \frac{1}{1 -2\tau +\tau^2} =
  \frac{1}{\operatorname{H}_R(-\tau)};
\end{equation*}
hence $R$ is Koszul by a result of L\"ofwall \thmcite[1.2]{CLf86}. If
$e$ is at least $3$, then $R$ is not Gorenstein and, therefore, $R$ is
Koszul by \thmcite[A]{LWCOVl07}. It is known that Koszul algebras are
quadratic, and the goal of this section is to prove that if $\k$ is
infinite, then a generic quadratic standard graded $\k$-algebra with
Hilbert series $1 + e\tau + (e-1)\tau^2$ has an exact zero divisor.

Recall that $R$ being quadratic means it is isomorphic to
$\poly{x_1,\ldots,x_e}/\mathfrak{q}$, where $\mathfrak{q}$ is an ideal
generated by homogeneous quadratic forms. The ideal $\mathfrak{q}$
corresponds to a subspace $V$ of the $\k$-vector space $W$ spanned by
$\setof{x_ix_j}{ 1 \le i \le j \le e}$. The dimension of $W$ is
$m=\frac{e(e+1)}{2}$, and since $R$ has Hilbert series $1 + e\tau +
(e-1)\tau^2$, the ideal $\mathfrak{q}$ is minimally generated by
$n=\frac{e^2-e+2}{2}$ quadratic forms; that is, $\dim[\k]{V}=n$.  In
this way, $R$ corresponds to a point in the Grassmannian
$\operatorname{Grass}_\k(n,m)$.

\begin{dfn}
  Let $e \ge 2$ be an integer and set $\mathcal{G}_\k(e) =
  \operatorname{Grass}_\k(n,m)$, where $n=\frac{e^2-e+2}{2}$ and $m =
  \frac{e(e+1)}{2}$. Points in $\mathcal{G}_\k(e)$ are in bijective
  correspondence with $\k$-algebras of embedding dimension $e$ whose
  defining ideal is minimally generated by $n$ homogeneous quadratic
  forms. For a point $\pi \in \mathcal{G}_\k(e)$ let $R_\pi$ denote
  the corresponding $\k$-algebra and let $\mathfrak{M}_\pi$ denote the
  irrelevant maximal ideal of $R_\pi$. Notice that
  $\operatorname{H}_{R_\pi}(\tau)$ has the form $1 +e\tau +
  (e-1)\tau^2+ \sum_{i=3}^\infty h_i\tau^i$ for every
  $\pi\in\mathcal{G}_\k(e)$.

  Consider the sets
  \begin{align*}
    \mathcal{E}_\k(e) &= \setof{\pi\in\mathcal{G}_\k(e)}{\text{there
        is an exact zero divisor in $R_\pi$}}\qand\\
    \mathcal{H}_\k(e) &=
    \setof{\pi\in\mathcal{G}_\k(e)}{\operatorname{H}_{R_\pi}(\tau) = 1
      +e\tau + (e-1)\tau^2}
  \end{align*}
  and recall that a subset of $\mathcal{G}_\k(e)$ is called open, if
  it maps to a Zariski open set under the Pl\"ucker embedding
  $\mathcal{G}_\k(e) \hookrightarrow \mathbb{P}_\k^{N}$, where $N =
  \binom{m}{n}-1$.
\end{dfn}

The sets $\mathcal{E}_\k(e)$ and $\mathcal{H}_\k(e)$ are non-empty:

\begin{exa}
  \label{exa:conca}
  Let $\k$ be a field and let $e \ge 2$ be an integer. The
  $\k$-algebra
  \begin{equation*}
    R = \frac{\poly{x_1,\ldots,x_e}}{(x_1^2) +
      (x_ix_j \mid 2 \le i \le j \le e)}
  \end{equation*}
  is local with Hilbert series $1 + e\tau + (e-1)\tau^2$. One has
  $\ann{x_1} = (x_1)$; in particular, $x_1$ is an exact zero divisor
  in $R$.
\end{exa}

\begin{thm}
  \label{thm:open}
  For every field $\k$ and every integer $e \ge 2$ the sets
  $\mathcal{E}_\k(e)$ and $\mathcal{H}_\k(e)$ are non-empty open
  subsets of the Grassmannian $\mathcal{G}_\k(e)$.
\end{thm}

\noindent The fact that $\mathcal{H}_\k(e)$ is open and non-empty is a
special case of \thmcite[1]{MHcDLk87}; for convenience a proof is
included below.

Recall that a property is said to hold for a generic algebra over an
infinite field if there is a non-empty open subset of an appropriate
Grassmannian such that every point in that subset corresponds to an
algebra that has the property.

\begin{cor}
  \label{cor:generic}
  Let $\k$ be an infinite field and let $e\ge 2$ be an integer. A
  generic standard graded $\k$-algebra with Hilbert series $1 + e\tau
  + (e-1)\tau^2$ has an exact zero~divisor.
\end{cor}

\begin{prf*}
  The assertion follows as the set $\mathcal{E}_\k(e) \cap
  \mathcal{H}_\k(e)$ is open by \thmref{open} and non-empty by
  \exaref{conca}.
\end{prf*}

\begin{rmk}
  Let $\Rm$ be a local ring. Following Avramov, Iyengar, and \c{S}ega
  (2008) we call an element $x$ in $R$ with $x^2=0$ and $x\m = \m^2$ a
  \emph{Conca generator}~of~$\m$. Conca proves in \seccite[4]{ACn00}
  that if $\k$ is algebraically closed, then the set
  \begin{equation*}
    \mathcal{C}_\k(e) =
    \setof{\pi\in\mathcal{G}_\k(e)}{\text{there is a Conca generator
        of $\mathfrak{M}_\pi$}}
  \end{equation*}
  is open and non-empty in $\mathcal{G}_\k(e)$.

  If $\m^3$ is zero, then it follows from \partlemref{ann}{c} that an
  element $x$ is a Conca generator of $\m$ if and only if the equality
  $\ann{x} = (x)$ holds. In particular, there is an inclusion
  \begin{equation*}
    \mathcal{C}_\k(e) \cap \mathcal{H}_\k(e) \subseteq
    \mathcal{E}_\k(e) \cap \mathcal{H}_\k(e).
  \end{equation*}
  If $\k$ is algebraically closed, then it follows from Conca's result
  combined with \thmref{open} and \exaref{conca} that both sets are
  non-empty and open in $\mathcal{G}_\k(e)$; in the next section we
  show that the inclusion may be strict.
\end{rmk}

\begin{rmk}
  \label{rmk:generic}
  Let $\Rmk$ be a local $\k$-algebra with $\m^3=0$ and assume that it
  is not Gorenstein. If $R$ admits a non-free totally reflexive
  module, then it has embedding dimension $e \ge 3$ and Hilbert series
  $1 + e\tau + (e-1)\tau^2$; see \pgref{m3tr}. Set
  \begin{equation*}
    \mathcal{T}_\k(e) =
    \setof{\pi\in\mathcal{H}_\k(e)}{\text{$R_\pi$ admits a non-free
        totally reflexive module}}.
  \end{equation*}
  We do not know if this is an open subset of $\mathcal{G}_\k(e)$, but
  it contains the non-empty open set $\mathcal{E}_\k(e) \cap
  \mathcal{H}_\k(e)$; hence the assertion in \thmref{generic}.
\end{rmk}

\begin{prf*}[Proof of \pgref{thm:open}]
  Let $\k$ be a field, let $e \ge 2$ be an integer, and let $S$ denote
  the standard graded polynomial algebra $\poly{x_1,\ldots,x_e}$. Set
  $N = \binom{m}{n} - 1$, where $m = \frac{e(e+1)}{2}$ and
  $n=\frac{e^2-e+2}{2}$. The Pl\"ucker embedding maps a point $\pi$ in
  the Grassmannian $\mathcal{G}_\k(e)$ to the point $(\mu_0(\pi):
  \cdots:\mu_N(\pi))$ in the projective space $\mathbb{P}_\k^N$, where
  $\mu_0(\pi), \ldots,\mu_N(\pi)$ are the maximal minors of any $m
  \times n$ matrix $\Pi$ corresponding to $\pi$. The columns of such a
  matrix give the coordinates, in the lexicographically ordered basis
  $\mathcal{B} = \setof{x_ix_j}{ 1 \le i \le j \le e}$ for $S_2$, of
  homogeneous quadratic forms $q_1,\ldots,q_n$. The algebra $R_\pi$ is
  the quotient ring $S/(q_1,\ldots,q_n)$.

  The sets $\mathcal{E}_\k(e)$ and $\mathcal{H}_\k(e)$ are non-empty
  by \exaref{conca}. Let $\Zeta$ be an $m \times n$ matrix of
  indeterminates and let $\chi_0,\ldots,\chi_N$ denote the maximal
  minors of $\Zeta$. We prove openness of each set $\mathcal{E}_\k(e)$
  and $\mathcal{H}_\k(e)$ in $\mathcal{G}_\k(e)$ by proving that the
  Pl\"ucker embedding maps it to the complement in $\mathbb{P}_\k^N$
  of the vanishing set for a finite collection of homogeneous
  polynomials in $\poly{\chi_0,\ldots,\chi_N}$.

  \subparagraph\emph{Openness of $\mathcal{E}_\k(e)$.}  Let $\pi$ be a
  point in $\mathcal{G}_\k(e)$, let $\Pi$ be a corresponding matrix,
  and let $\mathfrak{q}$ be the defining ideal for $R=R_\pi$. For a
  linear form $\ell = a_1x_1+\cdots+a_ex_e$ in $S$, let $l$ denote the
  image of $\ell$ in $R$. For $i\in\set{1,\ldots,e}$ let $[x_i\ell]$
  denote the column that gives the coordinates of $x_i\ell$ in the
  basis $\mathcal{B}$. Multiplication by $l$ defines a $\k$-linear map
  from $R_1$ to $R_2$. By assumption, one has $\dim[\k]{R_2} =
  \dim[\k]{R_1} - 1$, so the equality $lR_1 = R_2$ holds if and only
  if $l$ is annihilated by a unique, up to scalar multiplication,
  homogeneous linear form in $R$.  Set
  \begin{equation*}
    \Xi_\ell = \big(\, [x_1\ell] \mid [x_2\ell] \mid \cdots \mid [x_e\ell]
    \mid \Pi\,\big);
  \end{equation*}
  it is an $m \times (m+1)$ matrix with entries in $\k$. The equality
  $lR_1 = R_2$ holds if and only if the equality $\ell S_1 +
  \mathfrak{q}_2 = S_2$ holds, and the latter holds if and only if the
  matrix $\Xi_\ell$ has maximal rank. For $i\in\set{1,\ldots,m+1}$ let
  $\nu_i(\ell)$ denote the maximal minor obtained by omitting the
  $i$th column of $\Xi_\ell$. The columns of $\Pi$ are linearly
  independent, so $\Xi_\ell$ has maximal rank if and only if one of
  the minors $\nu_1(\ell),\ldots,\nu_e(\ell)$ is non-zero. Notice that
  each of the minors $\nu_1(\ell),\ldots,\nu_e(\ell)$ is a polynomial
  expression of degree $e-1$ in the coefficients $a_1,\ldots,a_e$ of
  $\ell$ and linear in the Pl\"ucker coordinates $\mu_1(\pi), \ldots,
  \mu_N(\pi)$. The column vector
  $(\nu_1(\ell)\;-\!\nu_2(\ell)\;\cdots\;(-1)^{m}\nu_{m+1}(\ell))^\mathrm{T}$
  is in the null-space of the matrix $\Xi_\ell$, so the element
  \begin{equation*}
    \bigg[ \Big( \sum_{i=1}^{e} (-1)^{i-1}\nu_i(\ell)x_i
    \Big)\ell\bigg] = \sum_{i=1}^e (-1)^{i-1}\nu_i(\ell)[x_i\ell]
  \end{equation*}
  is in the column space of the matrix $\Pi$. Set $\ell' =
  \sum_{i=1}^{e} (-1)^{i-1}\nu_i(\ell)x_i$; it follows that
  $\ell'\ell$ belongs to the ideal $\mathfrak{q}$, so one has $l'l=0$
  in $R$. By the discussion above, $l$ is the unique, up to scalar
  multiplication, annihilator in $R_1$ of $l'$ if and only if one of
  the maximal minors $\nu_1(\ell'),\ldots,\nu_e(\ell')$ of
  $\Xi_{\ell'}$ is non-zero. Moreover, if one of
  $\nu_1(\ell'),\ldots,\nu_e(\ell')$ is non-zero, then also one of the
  minors $\nu_1(\ell),\ldots,\nu_e(\ell)$ is non-zero by the
  definition of $\ell'$. Thus $l$ is an exact zero divisor in $R$ if
  and only if one of $\nu_1(\ell'),\ldots,\nu_e(\ell')$ is non-zero.
  Notice that each of these minors is a polynomial expression of
  degree $e$ in $\mu_1(\pi), \ldots, \mu_N(\pi)$ and degree
  $(e-1)^2$~in~$a_1,\ldots,a_e$.

  Let $\zeta_1, \ldots, \zeta_e$ be indeterminates and set $L =
  \zeta_1x_1 + \cdots + \zeta_ex_e$. For $i\in \set{1,\ldots,e}$ let
  $\nu_i$ denote the maximal minor of the matrix
  \begin{equation*}
    \big(\, [x_1L] \mid [x_2L] \mid \cdots \mid [x_eL]
    \mid \Zeta\,\big),
  \end{equation*}
  obtained by omitting the $i$th column. Each minor $\nu_i$ is a
  polynomial of degree $e-1$ in the indeterminates $\zeta_1, \ldots,
  \zeta_e$ and linear in $\chi_0,\ldots,\chi_N$. For $i\in
  \set{1,\ldots,e}$ set
  \begin{equation*}
    F_i = \nu_i(\nu_1, -\nu_2,\ldots,(-1)^{e-1}\nu_e);
  \end{equation*}
  each $F_i$ is a polynomial of degree $(e-1)^2$ in $\zeta_1, \ldots,
  \zeta_e$ and of degree $e$ in $\chi_0,\ldots,\chi_N$.  Consider
  $F_1,\ldots,F_e$ as polynomials in $\zeta_1, \ldots, \zeta_e$ with
  coefficients in $\poly{\chi_0,\ldots,\chi_N}$, and let $\mathcal{P}$
  denote the collection of these coefficients. The algebra $R$ has an
  exact zero divisor, i.e.\ the point $\pi$ belongs to
  $\mathcal{E}_\k(e)$, if and only if one of the polynomials $F_i$ in
  the algebra $(\poly{\chi_0,\ldots,\chi_N})[\zeta_1,\ldots,\zeta_e]$
  is non-zero, that is, if and only if the Pl\"ucker embedding maps
  $\pi$ to a point in the complement of the algebraic variety
  \begin{equation*}
    \bigcap_{P \in \mathcal{P}} \mathcal{Z}(P) \subseteq \mathbb{P}_\k^N.
  \end{equation*}

  \subparagraph \emph{Openness of $\mathcal{H}_\k(e)$.} Let $\pi$ be a
  point in $\mathcal{G}_\k(e)$, let $\Pi$ be a corresponding matrix,
  and let $\mathfrak{q} = (q_1,\ldots,q_n)$ be the defining ideal for
  $R_\pi$. Clearly, $\pi$ belongs to the subset $\mathcal{H}_\k(e)$ if
  and only if the equality $S_1\mathfrak{q}_2 = S_3$ holds.  Set $c =
  \binom{e+2}{3}$ and take as $\k$-basis for $S_3$ the $c$ homogeneous
  cubic monomials ordered lexicographically. For a homogeneous cubic
  form $f \in S_3$, let $[f]$ denote the column that gives its
  coordinates in this basis. The equality $S_1\mathfrak{q}_2 = S_3$
  holds if and only if the $c \times ne$ matrix
  \begin{equation*}
    \Xi = \big(\, [x_1q_1] \mid [x_1q_2] \mid \cdots \mid [x_1q_n] \mid [x_2q_1]
    \mid \cdots \mid [x_eq_n]\,\big)
  \end{equation*}
  has maximal rank, i.e.\ rank $c$ as one has $c \le ne$. Set
  \begin{equation*}
    \Epsilon = \big(\, [x_1(x_1^2)] \mid [x_1(x_1x_2)] \mid \cdots
    \mid [x_1(x_e^2)] \mid [x_2(x_1^2)] \mid \cdots \mid
    [x_e(x_e^2)]\,\big);
  \end{equation*}
  it is a $c \times me$ matrix, and each column of $\Epsilon$ is
  identical to a column in the $c\times c$ identity matrix
  $\Iota_c$. In particular, $\Epsilon$ has entries from the set
  $\set{0,1}$.  Let $\Delta$ be the matrix $\Pi^{\oplus e}$, that is,
  the block matrix with $e$ copies of $\Pi$ on the diagonal and $0$
  elsewhere; it is a matrix of size $me \times ne$. It is
  straightforward to verify the equality $\Xi = \Epsilon\Delta$. Set
  $g=ne-c$ and let $C$ denote the collection of all subsets of
  $\set{1,\ldots,me}$ that have cardinality $g$. For $i \in
  \set{1,\ldots,me}$ let $r_i$ denote the $i$th row of the identity
  matrix $\Iota_{me}$. For each $\gamma = \set{t_1,\ldots,t_g}$ in $C$
  set
  \begin{equation*}
    \Epsilon_\gamma = \left(
      \begin{array}{c}
        r_{t_1}\\ \vdots \\ r_{t_g}\\ \hline E
      \end{array} \right);
  \end{equation*}
  it is an $ne \times me$ matrix with entries from the set
  $\set{0,1}$.

  \subparagraph \emph{Claim.}  The matrix $\Xi$ has maximal rank if
  and only if there exists a $\gamma \in C$ such that
  $\Epsilon_\gamma\Delta$ has non-zero determinant.

  \emph{Proof.} Assume that the $ne\times ne$ matrix
  $\Epsilon_\gamma\Delta$ has non-zero determinant. One can write the
  determinant as a linear combination with coefficients in
  $\set{-1,0,1}$ of the $c$-minors of the submatrix $\Epsilon\Delta$,
  so $\Xi = \Epsilon\Delta$ has maximal rank. To prove the converse,
  assume that $\Epsilon\Delta$ has maximal rank, i.e.\ rank $c$. The
  matrix $\Delta$ has maximal rank, $ne$, so the rows of $\Delta$ span
  $\k^{ne}$. Therefore, one can choose $\gamma = \set{t_1,\ldots,t_g}$
  in $C$ such that rows number $t_1,\ldots,t_g$ in $\Delta$ together
  with the rows of $\Epsilon\Delta$ span $\k^{ne}$. That is, the rows
  of $\Epsilon_\gamma\Delta$ span $\k^{ne}$, so the determinant of
  $\Epsilon_\gamma\Delta$ is non-zero.

  \subparagraph The determinants $P_\gamma =
  \det{(\Epsilon_\gamma\Zeta^{\oplus e})}$, for $\gamma \in C$, yield
  $\binom{me}{g}$ homogeneous polynomials of degree $e$ in
  $\poly{\chi_0,\ldots,\chi_N}$. Under the Pl\"ucker embedding, $\pi$
  is mapped to a point in the complement of the algebraic variety
  \begin{equation*}
    \bigcap_{\gamma\in C} \mathcal{Z}(P_\gamma) \subseteq \mathbb{P}_\k^N
  \end{equation*}
  if and only if $\det{(\Epsilon_\gamma\Pi^{\oplus e})}$ is non-zero
  for some $\gamma \in C$. By Claim such a $\gamma$ exists if and only
  if $\Xi$ has maximal rank, that is, if and only if $\pi$ belongs to
  $\mathcal{H}_\k(e)$.
\end{prf*}

\begin{rmk}
  \label{rmk:zar}
  Let $\Rmk$ be a local $\k$-algebra with $\m^3=0$. The argument that
  shows the openness of $\mathcal{E}_\k(e)$ yields additional
  information. Namely, if there is an exact zero divisor in $R$, then
  one of the polynomials $F_i$ in the variables $\zeta_1, \ldots,
  \zeta_e$ is non-zero, and every point in the complement of its
  vanishing set $\mathcal{Z}(F_i) \subseteq \mathbb{P}_{\k}^{e-1}$
  corresponds to an exact zero divisor. Thus if $\k$ is infinite, then
  a generic homogeneous linear form in $R$ is an exact zero divisor.
\end{rmk}

\section{Short local rings without exact zero divisors---An example} %
\label{sec:example}

\noindent
Let $\k$ be a field; set
\begin{equation*}
  R=\poly[\k]{s,t,u,v}/(s^2,sv,t^2,tv,u^2,uv,v^2-st-su) \qand \m=(s,t,u,v)R.
\end{equation*}
Conca mentions in \cite[Example 12]{ACn00} that although $R$ is a
standard graded $\k$-algebra with Hilbert series $1 + 4\tau + 3\tau^2$
and $\ann{\m} = \m^2$, empirical evidence suggests that there is no
element $x\in R$ with $\ann{x} = (x)$; that is, there is no Conca
generator of~$\m$. \prpref{exa1} confirms this, and together with
\prpref{exa2} it exhibits properties of $R$ that frame the results in
the previous sections. In particular, these propositions show that
non-free totally reflexive modules may exist even in the absence of
exact zero divisors, and that exact zero divisors may exist also in
the absence of Conca generators.

\begin{prp}
  \label{prp:exa1}
  The following hold for the $\k$-algebra $R$.
  \begin{prt}
  \item There is no element $x$ in $R$ with $\ann{x} = (x)$.
  \item If $\,\k$ does not have characteristic $2$ or $3$, then the
    elements $s+t+2u-v$ and $3s+t-2u+4v$ form an exact pair of zero
    divisors in $R$.
  \item Assume that $\k$ has characteristic $3$.  If $\vartheta \in
    \k$ is not an element of the prime subfield $\mathsf{F}_3$, then
    the element $(1-\vartheta)s + \vartheta t + u + v$ is an exact
    zero divisor in~$R$. If $\,\k$ is $\mathsf{F}_3$, then there are
    no exact zero divisors in $R$.
  \item If $\,\k$ has characteristic $2$, then there are no exact zero
    divisors in $R$.
  \end{prt}
\end{prp}

\begin{prp}
  \label{prp:exa2}
  The $R$-module presented by the matrix
  \begin{equation*}
    \Phi=\begin{pmatrix}
      t&  -t+u-v\\
      t+u-v & s+u
    \end{pmatrix}
  \end{equation*}
  is indecomposable and totally reflexive, its first syzygy is
  presented by
  \begin{equation*}
    \Psi= \begin{pmatrix} -t+v & 2s+t-u+2v\\ t+u & s-u+v \end{pmatrix},
  \end{equation*}
  and its minimal free resolution is periodic of period 2.
\end{prp}

\begin{prf*}[Proof of \pgref{prp:exa1}]
  For an element $x = \alpha s + \beta t + \gamma u + \delta v$ in
  $\m$ we denote the images of $\alpha$, $\beta$, $\gamma$, and
  $\delta$ in $R/\m \is \k$ by $a$, $b$, $c$, and $d$.  For $x$ and
  $x' = \alpha' s + \beta' t + \gamma' u + \delta' v$ the product
  $xx'$ can be written in terms of the basis $\set{st,su,tu}$ for
  $\m^2$ as follows:
  \begin{equation}
    \label{eq:rr'}
    xx' = (ab'+ba' + dd')st + (ac' + ca' + dd')su + (bc' + cb')tu.
  \end{equation}
  The product $x\m$ is generated by the elements $xs = bst + csu$, $xt
  = ast + ctu$, $xu = asu + btu$, and $xv = dst + dsu$. By
  \rmkref{rank} the equality $x\m = \m^2$ holds if and only if the
  matrix
  \begin{equation*}
    \Xi_x=
    \begin{pmatrix}
      b & a & 0 & d\\ c & 0 & a & d\\ 0 & c & b &0
    \end{pmatrix}
  \end{equation*}
  has a non-zero $3$-minor; that is, if and only if one of
  \begin{equation}
    \label{eq:minors}
    \begin{split}
      \begin{aligned}
        \mu_1(x) &= -2abc, & \mu_2(x)&=cd(c-b),\\
        \mu_3(x) & =bd(c-b), \qand& \mu_4(x)&=-ad(c+b)
      \end{aligned}
    \end{split}
  \end{equation}
  is non-zero.  Thus elements $x$ and $x'$ with $xx'=0$ form an exact
  pair of zero divisors if and only if one of the minors $\mu_1(x),
  \ldots, \mu_4(x)$ is non-zero, and one of the minors $\mu_1(x'),
  \ldots, \mu_4(x')$ of the matrix $\Xi_{x'}$ is non-zero.

  \subparagraph (a): For the equality $\ann{x} = (x)$ to hold, $x$
  must be an element in $\m\backslash\m^2$. If $x = \alpha s + \beta t
  + \gamma u + \delta v$ satisfies $x^2=0$, then \eqref{rr'} yields
  \begin{equation*}
    2ab+d^2=0, \quad 2ac+d^2=0, \qand 2bc=0.
  \end{equation*}
  It follows that $b$ or $c$ is zero and then that $d$ is zero. Thus
  each of the minors $\mu_1(x), \ldots, \mu_4(x)$ is zero, whence $x$
  does not generate $\ann{x}$.

  \subparagraph (b): For the elements $x=s+t+2u-v$ and $x' =
  3s+t-2u+4v$, it is immediate from \eqref{rr'} that the product $xx'$
  is zero, while \eqref{minors} yields $\mu_3(x) =-1$ and
  $\mu_1(x')=12$. If $\k$ does not have characteristic $2$ or $3$,
  then $\mu_1(x')$ is non-zero, so $x$ and $x'$ form an exact pair of
  zero divisors in $R$.

  \subparagraph (c): Assume that $\k$ has characteristic $3$.  If $\k$
  properly contains $\mathsf{F}_3$, then choose an element $\vartheta
  \in\k\backslash \mathsf{F}_3$ and set
  \begin{equation*}
    x = (1-\vartheta)s + \vartheta t + u + v.
  \end{equation*}
  Observe that $\mu_2(x) = 1-\vartheta$ is non-zero. If $\vartheta$ is
  not a 4th root of unity, set
  \begin{equation*}
    x'=(1 + \vartheta)s + t - \vartheta u -(1+\vartheta^2)v,
  \end{equation*}
  and note that the minor $\mu_2(x') =
  \vartheta(1+\vartheta)(1+\vartheta^2)$ is non-zero. From \eqref{rr'}
  one readily gets $xx'=0$, so $x$ and $x'$ form an exact pair of zero
  divisors.  If $\vartheta$ is a 4th root of unity, then one has
  $\vartheta^2 = -1$. Set
  \begin{equation*}
    x'' = (1-\vartheta)s + t + \vartheta u - v;
  \end{equation*}
  then the minor $\mu_2(x'') = \vartheta(1-\vartheta)$ is non-zero,
  and it is again straightforward to verify the equality
  $xx''=0$. Therefore, $x$ and $x''$ form an exact pair of zero
  divisors.

  Assume now $\k=\mathsf{F_3}$ and assume that the elements $x =
  \alpha s + \beta t + \gamma u + \delta v$ and $x' = \alpha' s +
  \beta' t + \gamma' u + \delta' v$ form an exact pair of zero
  divisors in $R$. One of the minors $\mu_1(x), \ldots, \mu_4(x)$ is
  non-zero, and one of the minors $\mu_1(x'), \ldots, \mu_4(x')$ is
  non-zero, so it follows from \eqref{minors} that $abc$ or $d$ is
  non-zero and that $a'b'c'$ or $d'$ is non-zero.

  First assume $dd'\ne 0$; we will show that the six elements
  $a,a',b,b',c,c'$ are non-zero and derive a contradiction. Suppose
  $b=0$, then \eqref{rr'} yields $cb'=0$ and $ab' \ne 0$, which forces
  $c=0$. This, however, contradicts the assumption that one of the
  minors $\mu_1(x), \ldots, \mu_4(x)$ is non-zero, cf.~\eqref{minors}.
  Therefore, $b$ is non-zero. A parallel arguments show that $c$ is
  non-zero, and by symmetry the elements $b'$ and $c'$ are
  non-zero. Suppose $a=0$, then it follows from \eqref{rr'} that $a'$
  is non-zero, as $dd'\ne 0$ by assumption. However, \eqref{rr'} also
  yields
  \begin{equation*}
    a(b'-c') = a'(c-b),
  \end{equation*}
  so the assumption $a=0$ forces $c=b$, which contradicts the
  assumption that one of the minors $\mu_1(x), \ldots, \mu_4(x)$ is
  non-zero. Thus $a$ is non-zero, and by symmetry also $a'$ is
  non-zero. Without loss of generality, assume $a=1=a'$, then
  \eqref{rr'} yields
  \begin{equation}
    \label{eq:minor3}
    b'+ b = c'+c \qand bc'+cb' =0.
  \end{equation}
  Eliminate $b'$ between these two equalities to get
  \begin{equation}
    \label{eq:bccb}
    b(c'-c) + c(c'+c) = 0.
  \end{equation}
  As $c$ and $c'$ are non-zero elements in $\mathsf{F}_3$, the
  elements $c'-c$ and $c'+c$ are distinct, and their product is
  $0$. Thus one and only one of them is $0$, which contradicts
  \eqref{bccb} as both $b$ and $c$ are non-zero.

  Now assume $dd'=0$. Without loss of generality, assume that $d$ is
  zero, then $abc$ is non-zero. It follows from \eqref{minors} that
  the elements $a'$, $b'$, and $c'$ cannot all be zero, and then
  \eqref{rr'} shows that they are all non-zero. Thus all six elements
  $a,a',b,b',c,c'$ are non-zero, and as above this leads to a
  contradiction.

  \subparagraph (d): Assume that $\k$ has characteristic $2$ and that
  the elements $x = \alpha s + \beta t + \gamma u + \delta v$ and $x'
  = \alpha' s + \beta' t + \gamma' u + \delta' v$ form an exact pair
  of zero divisors in $R$. From \eqref{minors} one gets $d\ne 0$,
  $d'\ne 0$, $b\ne c$, and $b'\ne c'$.  Arguing as in part (c), it is
  straightforward to verify that the six elements $a,a',b,b',c,c'$ are
  non-zero. Without loss of generality, assume $a=1=a'$. From
  \eqref{rr'} the following equalities emerge:
  \begin{equation*}
    b'+b+dd' = 0,\quad c'+c+dd'=0,\qand bc'+cb' =0.
  \end{equation*}
  The first two equalities yield $b'+ b \ne 0$ and $c' +c \ne 0$.
  Further, elimination of $dd'$ yields $b' = b+c+c'$. Substitute this
  into the third equality to get
  \begin{equation*}
    b(c'+c) + c(c'+c) = 0.
  \end{equation*}
  As $c'+c$ is non-zero this implies $b=c$, which is a contradiction.
\end{prf*}

\begin{prf*}[Proof of \pgref{prp:exa2}]
  It is easy to verify that the products $\Phi\Psi$ and $\Psi\Phi$ are
  zero, whence
  \begin{equation*}
    F:\quad \cdots \lra R^2 \xra{\Phi} R^2 \xra{\Psi} R^2 \xra{\Phi}
    R^2 \lra \cdots,
  \end{equation*}
  is a complex.  We shall first prove that $F$ is totally acyclic. To
  see that it is acyclic, one must verify the equalities $\Im{\Psi} =
  \Ker{\Phi}$ and $\Im{\Phi} = \Ker{\Psi}$.  Assume that the element
  $\binom{x}{x'}$ is in $\Ker\Phi$. It is straightforward to verify
  that $\m^2R^2$ is contained in the image of $\Psi$, so we may assume
  that $x$ and $x'$ have the form $x=as+bt+cu+dv$ and
  $x'=a's+b't+c'u+d'v$, where $a,b,c,d$ and $a',b',c',d'$ are elements
  in $\k$.  Using that $\set{st,su,tu}$ is a basis for $\m^2$, the
  assumption $\Phi\binom{x}{x'}=0$ can be translated into the
  following system of equations
  \begin{equation*}
    0=
    \begin{cases}
      a - a' - d'\\
      a' - d'\\
      c + b' -c'\\
      a - d +b'\\
      a-d+a'+c'\\
      b+c+b'
    \end{cases}
  \end{equation*}
  From here one derives, in order, the following identities
  \begin{equation*}
    a'=d', \ \ a=2a', \ \ d=2a'+b', \ \ c'=-a'+b', \
    c=-a', \ \text{ and } \ b=a'-b',
  \end{equation*}
  which immediately yield $\binom{x}{x'} = \Psi\binom{b'}{a'}$. This
  proves the equality $\Im{\Psi}=\Ker{\Phi}$. Similarly, it is easy to
  check that $\m^2R^2$ is contained in $\Im{\Phi}$, and for an element
  \begin{equation}
    \label{eq:element}
    \binom{x}{x'} = \binom{as+bt+cu+dv}{a's+b't+c'u+d'v},
  \end{equation}
  where $a,a',\ldots,d,d'$ are elements in $\k$, one finds that
  $\Psi\binom{x}{x'}=0$ implies $\binom{x}{x'} = \Phi\binom{b'}{a'}$.
  This proves the equality $\Im{\Phi} = \Ker{\Psi}$, so $F$ is
  acyclic. The differentials in the dual complex $\Hom{F}{R}$ are
  represented by the matrices $\Phi^\mathrm{T}$ and
  $\Psi^\mathrm{T}$. One easily checks the inclusions $\m^2R^2
  \subseteq \Im{\Phi^\mathrm{T}}$ and $\m^2R^2 \subseteq
  \Im{\Psi^\mathrm{T}}$. Moreover, for an element of the form
  \eqref{element} one finds
  \begin{alignat*}{3}\textstyle
    \Phi^\mathrm{T}\binom{x}{x'} &=0 & &\quad\text{implies}\quad&
    \textstyle\binom{x}{x'} &=
    \textstyle\Psi^\mathrm{T}\binom{b'}{a'-2b'}\qand\\
    \textstyle\Psi^\mathrm{T}\binom{x}{x'} &=0 &
    &\quad\text{implies}\quad& \textstyle\binom{x}{x'} &=
    \textstyle\Phi^\mathrm{T}\binom{-b'}{a'}.
  \end{alignat*}
  This proves that also $\Hom{F}{R}$ is acyclic, so the module, $M$,
  presented by $\Phi$ is totally reflexive. Moreover, the first syzygy
  of $M$ is presented by $\Psi$, and the minimal free resolution of
  $M$ is periodic of period $2$.

  To prove that $M$ is indecomposable, assume that there exist
  matrices $\Alpha=(a_{ij})$ and $\Beta=(b_{ij})$ in $\GL[R]{2}$, such
  that $\Alpha\Phi\Beta$ is a diagonal matrix; i.e.\ the equalities
  \begin{align}
    \label{eq:dia1}
    \begin{split}
      0 &= (a_{12}b_{22})s + ((a_{11}+a_{12})b_{12} - a_{11}b_{22})t +
      (a_{12}b_{12} + (a_{11}+a_{12})b_{22})u\\&\quad - (a_{12}b_{12}
      + a_{11}b_{22})v
    \end{split}
  \end{align}
  and
  \begin{align}
    \label{eq:dia2}
    \begin{split}
      0 &= (a_{22}b_{21})s + ((a_{21}+a_{22})b_{11} - a_{21}b_{21})t
      +(a_{22}b_{11} + (a_{21}+a_{22})b_{21})u\\&\quad -
      (a_{22}b_{11}+a_{21}b_{21})v
    \end{split}
  \end{align}
  hold.  As the matrices $\Alpha$ and $\Beta$ are invertible, neither
  has a row or a column with both entries in $\m$. Since the elements
  $s$, $t$, $u$, and $v$ are linearly independent modulo $\m^2$, it
  follows from \eqref{dia1} that $a_{12}b_{22}$ is in $\m$. Assume
  that $a_{12}$ is in $\m$, then $a_{11}$ and $a_{22}$ are not in
  $\m$. It also follows from \eqref{dia1} that $a_{12}b_{12} +
  a_{11}b_{22}$ is in $\m$, which forces the conclusion
  $b_{22}\in\m$. However, this implies that $b_{21}$ is not in $\m$,
  so $a_{22}b_{21}$ is not in $\m$ which contradicts \eqref{dia2}. A
  parallel argument shows that also the assumption $b_{22}\in\m$ leads
  to a contradiction. Thus $M$ is indecomposable.
\end{prf*}

\section{Families of non-isomorphic modules of infinite length}
\label{sec:ext}

\noindent
Most available proofs of the existence of infinite families of totally
reflexive modules are non-constructive. In the previous sections we
have presented constructions that apply to local rings with exact zero
divisors. In \cite{HHl11} Holm gives a different construction; it
applies to rings of positive dimension which have a special kind of
exact zero divisors. Here we provide one that does not depend on exact
zero~divisors.

\begin{con}
  \label{con:po}
  Let $\Rm$ be a local ring and let $\varkappa = \set{x_1,\dots,x_e}$
  be a minimal set of generators for $\m$. Let $N$ be a finitely
  generated $R$-module and let $N_1$ be its first syzygy. Let $F \onto
  N$ be a projective cover, and consider the element
  \begin{align*}
    \xi:\quad 0 \lra N_1 \xra{\iota} F \lra N \lra 0
  \end{align*}
  in $\Ext{1}{N}{N_1}$. For $i\in\set{1,\ldots,e}$ and $j\in \NN$
  recall that $x_i^j\xi$ is the second row in the diagram
  \newlength{\dimb}\settowidth{\dimb}{$x_i^j$}
  \begin{equation*}
    \xymatrix{
      \hspace{\dimb}\xi\colon\mspace{-50mu} & 0
      \ar[r]
      & N_1
      \ar[r]^\iota
      \ar[d]^{x_i^j}
      & F
      \ar[r]
      \ar[d]
      & N
      \ar[r]
      \ar@{=}[d]
      & 0
      \\
      x_i^j\xi\colon \mspace{-50mu}& 0
      \ar[r]
      & N_1
      \ar[r]^-{\omega^{(i,j)}}
      & \Po{j}
      \ar[r]
      & N
      \ar[r]
      & 0,
    }
  \end{equation*}
  where the left-hand square is the pushout of $\iota$ along the
  multiplication map $x_i^j$. The diagram defines $\Po{j}$ uniquely up
  to isomorphism of $R$-modules. Set
  \begin{equation*}
    \mathcal{P}(\varkappa\,;N) = \setof{\Po{j}}{1 \le i \le e,\ j\in\NN};
  \end{equation*}
  note that every module in $\mathcal{P}(\varkappa\,;N)$ can be
  generated by $\bet{0}{N} + \bet{1}{N}$ elements.
\end{con}

\begin{lem}
  \label{lem:po}
  Let $\Rm$ be a local ring and let $N$ be a finitely generated
  $R$-module. If, for some minimal set $\varkappa =
  \set{x_1,\dots,x_e}$ of generators for $\m$, the set
  $\mathcal{P}(\varkappa\,;N)$ contains only finitely many pairwise
  non-isomorphic modules, then the $R$-module $\Ext{1}{N}{N_1}$ has
  finite length.
\end{lem}

\begin{prf*}
  Let $\varkappa = \set{x_1,\dots,x_e}$ be a minimal set of generators
  for $\m$ and assume that $\mathcal{P}(\varkappa\,;N)$ contains only
  finitely many pairwise non-isomorphic modules. Given an index
  $i\in\set{1,\ldots,e}$, there exist positive integers $m$ and $n$
  such that the modules $\Po{m}$ and $\Po{n+m}$ are isomorphic. Since
  $x_i^{n+m}\xi$ equals $x_i^n(x_i^m\xi)$, the module $\Po{n+m}$ comes
  from the pushout of $\omega^{(i,m)}$ along the multiplication
  map~$x_i^n$; cf.~\conref{po}. Thus there is an exact sequence
  \begin{equation*}
    0 \lra N_1 \xra{\alpha} N_1 \oplus \Po{m} \lra \Po{n+m} \lra 0,
  \end{equation*}
  where $\alpha = \begin{pmatrix} x_i^n &\mspace{-8mu}
    -\omega^{(i,m)} \end{pmatrix}$. It follows from the isomorphism
  $\Po{m} \is \Po{n+m}$ and Miyata's theorem \cite{TMi67} that this
  sequence splits. Hence, it induces a split monomorphism
  \begin{equation*}
    \Ext{1}{N}{N_1} \xra{\;\Ext{1}{N}{\alpha}\;}
    \Ext{1}{N}{N_1} \oplus \Ext{1}{N}{\Po{m}}.
  \end{equation*}
  Let $\beta$ be a left-inverse of $\Ext{1}{N}{\alpha}$, set $m_i = m$
  and notice that the element
  \begin{equation*}
    x_i^{m_i}\xi = \beta\Ext{1}{N}{\alpha}(x_i^{m_i}\xi) = \beta
    \begin{pmatrix}
      x_i^{n+m_i}\xi &\mspace{-8mu} 0
    \end{pmatrix}
    = x_i^{n+m_i}\beta
    \begin{pmatrix}
      \xi &\mspace{-8mu} 0
    \end{pmatrix}
  \end{equation*}
  belongs to $\m^{n+m_i}\Ext{1}{N}{N_1}$.

  For every index $i\in\set{1,\ldots,e}$ let $m_i$ be the positive
  integer obtained above. With $h = m_1 + \cdots + m_e$ there is an
  inclusion $\m^h\Ext{1}{N}{N_1} \subseteq \m^{h+1}\Ext{1}{N}{N_1}$,
  so Nakayama's lemma yields $\m^h\Ext{1}{N}{N_1} =0$.
\end{prf*}

\begin{thm}
  Let $\Rm$ be a local ring and let $\varkappa = \set{x_1,\dots,x_e}$
  be a minimal set of generators for $\m$. If there exists a totally
  reflexive $R$-module $N$ and a prime ideal $\p \ne \m$ such that
  $N_\p$ is not free over $R_\p$, then the set
  $\mathcal{P}(\varkappa\,;N)$ contains infinitely many indecomposable
  and pairwise non-isomorphic totally reflexive $R$-modules.
\end{thm}

\begin{prf*}
  It follows from the assumptions on $N$ that every module in the set
  $\mathcal{P}(\varkappa\,;N)$ is totally reflexive and that the
  $R$-module $\Ext{1}{N}{N_1}$ has infinite length, as its support
  contains the prime ideal $\p\ne\m$. By \lemref[]{po} the set
  $\mathcal{P}(\varkappa\,;N)$ contains infinitely many pairwise
  non-isomorphic modules. Every module in $\mathcal{P}(\varkappa\,;N)$
  is minimally generated by at most $\bet{0}{N} + \bet{1}{N}$
  elements; see~\conref{po}. Therefore, every infinite collection of
  pairwise non-isomorphic modules in $\mathcal{P}(\varkappa\,;N)$
  contains infinitely many indecomposable modules.
\end{prf*}

\section*{Acknowledgments} 

\noindent
It is our pleasure to thank Manoj Kummini and Christopher Monico for
helpful conversations related to some of the material in this paper.


\bibliographystyle{amsplain} %

\providecommand{\arxiv}[2][AC]{\mbox{\href{http://arxiv.org/abs/#2}{\sf
      arXiv:#2 [math.#1]}}}
\providecommand{\oldarxiv}[2][AC]{\mbox{\href{http://arxiv.org/abs/math/#2}{%
      \sf arXiv:math/#2
      [math.#1]}}}\providecommand{\MR}[1]{\mbox{\href{http://www.ams.org/mathscine%
      t-getitem?mr=#1}{#1}}}
\renewcommand{\MR}[1]{\mbox{\href{http://www.ams.org/mathscinet-getitem?mr=#%
      1}{#1}}} \providecommand{\bysame}{\leavevmode\hbox
  to3em{\hrulefill}\thinspace}
\providecommand{\MR}{\relax\ifhmode\unskip\space\fi MR }
\providecommand{\MRhref}[2]{%
  \href{http://www.ams.org/mathscinet-getitem?mr=#1}{#2} }
\providecommand{\href}[2]{#2}

\end{document}